\documentclass{article}
\usepackage[utf8]{inputenc}

\usepackage[english]{babel}

\usepackage[utf8]{inputenc}
\usepackage[T1]{fontenc}
\usepackage[normalem]{ulem}

\usepackage{amssymb}
\usepackage{amsmath,amsfonts}
\usepackage[tikz]{bclogo}
\usetikzlibrary{arrows, arrows.meta}
\usepackage[top=3cm,bottom=3cm,left=2.5cm,right=2.5cm]{geometry}

\usepackage{lmodern}
\usepackage{textcomp}
\usepackage{mathtools, bm}
\usepackage{amssymb, bm}
\usepackage[unq]{unique}
\usepackage{enumerate}
\usepackage{comment}
\usepackage[breaklinks]{hyperref}

\usepackage[numbers]{natbib}  

\usepackage[breaklinks]{hyperref}
\usepackage[numbers]{natbib}
\usepackage{graphicx} % Required for inserting images

\usepackage{amsthm}

\usepackage{authblk}
\usepackage{tikz}
\usepackage{caption}
\usepackage{subcaption}

\usepackage{cleveref}

\newtheorem{theorem}{Theorem}[section]
\newtheorem{lemma}[theorem]{Lemma}
\newtheorem{proposition}[theorem]{Proposition}

\newtheorem{claim}[theorem]{Claim}

\theoremstyle{definition}
\newtheorem{remark}[theorem]{Remark}
\newtheorem{conjecture}[theorem]{Conjecture}

\newtheorem{example}[theorem]{Example}
\newtheorem{definition}[theorem]{Definition}

\newcommand{\EE}{\mathbb{E}}

\newcommand{\defined}{\mathrel{\coloneqq}}
\undef{\set}
\DeclarePairedDelimiter{\set}{\lbrace}{\rbrace}
\newcommand{\eps}{\varepsilon}
\newcommand{\sseq}{\subseteq}

\newcommand{\union}{\mathbin{\cup}}

\newcommand{\cS}{\mathcal{S}}

\newcommand{\bc}{\mathbf{c}}
\newcommand{\from}{\colon}
\undef{\abs}
\DeclarePairedDelimiterX{\abs}[1]
{\lvert}{\rvert}{\ifblank{#1}{\,\cdot\,}{#1}}

\DeclareMathOperator{\Var}{Var}

\title{Discrepancies of spanning trees in dense graphs}

\author[1]{Lawrence Hollom}
\author[2]{Lyuben Lichev}
\author[1]{Adva Mond}
\author[1]{Julien Portier}
\affil[1]{Department of Pure Mathematics and Mathematical Statistics, University of Cambridge, Cambridge, CB3 0WA, United Kingdom}
\affil[2]{Institute of Science and Technology Austria (ISTA), 3400 Klosterneuburg, Austria}

\date{\today}

\begin{document}

\maketitle

\begin{abstract}
%We first consider discrepancy problems in a setting introduced by Erd\H{o}s, F\"{u}redi, Loebl and S\'{o}s. 
We address several related problems on combinatorial discrepancy of trees in a setting introduced by Erd\H{o}s, F\"{u}redi, Loebl and S\'{o}s.
Given a fixed tree $T$ on $n$ vertices and an edge-colouring of the complete graph $K_n$, for every colour, we find a copy of $T$ in $K_n$ where the number of edges in that colour significantly exceeds its expected count in a uniformly random embedding.
In particular, this resolves a problem posed by Erd\H{o}s, F\"{u}redi, Loebl and S\'{o}s by generalising their work from two to many colours.
Furthermore, if $T$ has maximum degree $\Delta\leq\eps n$ for sufficiently small $\eps > 0$ and the edge-colouring of $K_n$ is both balanced and ``not too close'' to one particular instance, we show that, for every colour, there is a copy of $T$ in $K_n$ where that colour appears on linearly more edges than any other colour.
Several related examples are provided to demonstrate the necessity of the introduced structural restrictions.
Moreover, when $\Delta$ is a constant, we extend these results to sufficiently dense host graphs in place of $K_n$.
Our proofs combine saturation arguments for the existence of particular coloured substructures and analysis of conveniently defined local exchanges.
	
Using similar methods, we investigate the existence of copies of a graph $H$ with prescribed number of edges in each colour in $2$-edge-coloured dense host graphs.
In particular, for a graph $H$ with bounded maximum degree and balanced $2$-edge-colourings $\mathbf{c}$ of a host graph $G$ with minimum degree at least $(1-\eps)n$ for some $\eps > 0$, we show that, for any sufficiently large $n$ and sufficiently small $\eps$, there exists a copy of $H$ where the number of edges in the two colours differ by at most $2$.
Moreover, we completely characterise the pairs $(H,\mathbf{c})$ for which the difference of $2$ cannot be improved.
As a consequence, we refute a conjecture by Mohr, Pardey, and Rautenbach.
\end{abstract}

\section{Introduction}

Discrepancy theory is a relatively modern and very active subfield of mathematical optimisation and algorithm design.
Its goal is to study the existence and the properties of unbalanced substructures of some weighted or coloured mathematical structure.
The roots of the subject may be found in the pioneering work of Weyl~\cite{Wey16}, and relations to classic mathematical fields like number theory, combinatorics, ergodic theory and discrete geometry were exhibited shortly after.
For a broad introduction to the subject, we recommend the outstanding books of Beck and Chen~\cite{BC87}, Matou\v{s}ek~\cite{Mat99} and Chazelle~\cite{Cha00}.

The current work contributes to a line of research initiated by Erd\H{o}s, F\"uredi, Loebl and S\'os~\cite{EFLS95} who showed that there is a constant $c > 0$ such that, for every tree $T$ on $n$ vertices and maximum degree $\Delta$ and every 2-colouring (that is, colouring in two colours) of the edges of $K_n$, one can find a copy of $T$ where one of the colours is given to at least $c(n-1-\Delta)$ more edges than the other colour.
Alternatively, in the language of discrepancy, the above result claims that some copy of $T$ has discrepancy at least $c(n-1-\Delta)$.

In this two-colour setting, the discrepancy of the tree can be viewed as the `imbalance' between the two colours, that is, the difference of the sizes of their pre-images. However, there is more than one natural way to adjust this notion to more colours.
In this paper, we consider several of them.

Our results go in two different directions.
First, in any edge-colouring of $K_n$ in $r\ge 2$ colours, we generalise the Erd\H{o}s--F\"uredi--Loebl--S\'os theorem by exhibiting copies of a tree $T$ where the number of edges in some colour is significantly over-represented. 
Later, we extend our study to approximately balanced edge-colourings of dense host graphs ($K_n$ being a particular example) and find copies of $T$ where a predetermined colour class contains significantly more edges than any other colour class.
%the frequency a \emph{fixed} colour in optimised.
Along the way, we exhibit several curious pathological constructions which indicate that certain structural constraints we rely on are actually needed.
These results are described in Sections~\ref{subsubsec:large-imbalance-intro1} and~\ref{subsubsec:large-imbalance-intro}.

In a different direction, for a given graph $H$ with constant maximum degree, under certain constraints on a $2$-colouring of the host graph $G$, we exhibit copies of $H$ with a prescribed number of edges in each colour.
The relevant results are described in \Cref{subsubsec:k-sum-intro}.

\subsection{Maximising the number of edges in some colour}
\label{subsubsec:large-imbalance-intro1}

Our first result generalises the main theorem of Erd\H{o}s, F\"uredi, Loebl and S\'os~\cite{EFLS95} to $r \ge 2$ colours, thereby resolving a problem raised in the same paper.

\begin{theorem}
		\label{col:discErdos}
		There exists a universal constant $c>0$ such that the following holds. Fix an integer $r\ge 2$ and an $r$-colouring of the edges of $K_n$. 
		Then, for every tree $T$ on $n$ vertices with maximum degree $\Delta$, there exists a colour $i\in [r]$ and a copy of $T$ in which the number of edges in colour i is at least
        $$\frac{n-1}{r}+\frac{c}{r^2}(n-1-\Delta).$$
\end{theorem}

In fact, we derive \Cref{col:discErdos} from a more general result, where we generalise the main result from~\cite{EFLS95} in two ways: we consider colourings in $r\ge 2$ colours and, for \emph{each} of the colours involved, we find a copy of $T$ where this colour is over-represented compared to the average number of edges it is given to in a uniformly random embedding.
\Cref{col:discErdos} follows from the next result applied to some most commonly met colour.

\begin{theorem}\label{thm:col_disc_1}
		There exists a universal constant $C > 0$ such that the following holds.
		Fix an integer $r\ge 2$ and an $r$-colouring of the edges of $K_n$. 
		For every colour $i\in [r]$, denote by $m_i$ the number of edges in colour~$i$.
		Then, for every $i\in [r]$ and every tree $T$ on $n$ vertices with maximum degree $\Delta$,
		there is a copy of $T$ with at least
		\[\frac{2m_i}{n} + C(n-1-\Delta)\left(\frac{m^*_i}{n^2} \right)^2 \]
		edges in colour $i$, where $m^*_i \coloneqq \min\left(m_i, \binom{n}{2}-m_i \right)$.
\end{theorem}

\noindent
Note that Theorem~\ref{thm:col_disc_1} can be used to find copies of the tree $T$ where the colour $i$ is under-represented instead: indeed, one can group the colours in $[r]\setminus \{i\}$ together into a new super-colour and apply Theorem~\ref{thm:col_disc_1} for that super-colour.

The topic of discrepancy in edge-coloured graphs has received considerable attention in the last few years. 
Balogh, Csaba, Jing and Pluh\'ar~\cite{BCJP20} studied the discrepancy of spanning trees and Hamilton cycles in 2-colourings of graphs with high minimum degree.
Notably, they showed that every graph on $n$ vertices with minimum degree at least $(3/4+\eps)n$ with any $\eps\in (0,1/4)$ contains a Hamilton cycle with discrepancy at least $\eps n/32$.
This result was later generalised independently by Freschi, Hyde, Lada and Treglown~\cite{FHLT21} and Gishboliner, Krivelevich and Michaeli~\cite{GKM22b} to more colours.
In particular, for every $r\ge 2$, both works showed that every $r$-coloured graph on $n$ vertices with minimum degree at least $(1/2+1/2r)n + d$ contains a Hamilton cycle where some colour is given to a least $n/r+cd$ edges where $c=1/6r^2$ in~\cite{FHLT21} and $c=2$ in~\cite{GKM22b}.
Moreover,~\cite{GKM22b} further generalised the results on spanning trees from~\cite{BCJP20} by exhibiting weak connectivity conditions for the host graph which are sufficient to provide the existence of a spanning tree with high discrepancy for any number of colours (for a suitable extension of discrepancy to the multi-colour case).

Gishboliner, Krivelevich and Michaeli~\cite{GKM22a} also studied $r$-colourings of the Erd\H{o}s-R\'enyi random graph $G(n,p)$ above its sharp threshold for Hamiltonicity $p_{\mathrm{Ham}} = (\log n+\log\log n)/n$.
They proved that, with high probability, for every $r$-colouring of the edges of $G(n,p)$ with $n(p-p_{\mathrm{Ham}})\to \infty$, there is a Hamilton cycle with at least $(2/(r+1)-o(1))n$ edges in the same colour.
A similar result was also provided for perfect matchings.
In the 2-colour setting, the discrepancy of powers of Hamilton cycles was further examined by Brada\v{c}~\cite{Bra22}. 
Balogh, Csaba, Pluh\'ar and Treglown~\cite{BCPT21} studied the case of $K_r$-factors and obtained a discrepancy version of the Hajnal-Szemerédi theorem by identifying the minimum degree threshold for the existence of a $K_r$-factor with high discrepancy.
Brada\v{c}, Christoph and Gishboliner~\cite{BCG2023} then generalised their result and, for any graph $H$, found tight minimum degree conditions for the existence of an $H$-factor with linear discrepancy.
Krishna, Michaeli, Sarantis, Wang and Wang~\cite{KMSWW23} generalised graph discrepancy in the 2-colour setting by introducing high-dimensional discrepancy where every edge is associated to a unit $d$-dimensional vector for some $d\ge 1$.
They investigated high-dimensional discrepancy in a setting where the host graph is an $n$-vertex tree $T$ and the subgraphs of interest are all subtrees of $T$.
Finally, the discrepancy of tight Hamilton cycles in $k$-uniform hypergraphs was the main focus of a recent work by Gishboliner, Glock and Sgueglia~\cite{GGS23}.

\subsection{Maximising the number of edges in a given colour}
\label{subsubsec:large-imbalance-intro}

While Theorem~\ref{thm:col_disc_1} guarantees the existence of a copy of $T$ where a certain colour is over-represented, it does not give any information about the proportion of edges in any other colour.
The following natural question arises: if the edges of $G$ are approximately evenly distributed among all colours, does there necessarily exist a copy of the tree $T$ where a fixed colour is met more frequently than any other colour? Our next result shows that this is not the case, even when restricted to trees with maximum degree $\Delta \leq \eps n$ for any small $\eps >0$.

We say than an $r$-colouring of a graph is \emph{balanced} if the number of edges in any two different colours differs by at most 1. 
Also, for an integer $\ell\ge 1$, an \emph{$\ell$-caterpillar} is a tree where every vertex lies within distance 1 of a central path of length $\ell$, and the degrees of the vertices in this path differ by at most 1.

\begin{proposition}\label{prop:example_many_colours}
	Fix $\eps\in (0,1)$ sufficiently small. Then, for all integers $r\ge 1/\eps^5$, the following holds:
	for all sufficiently large $n$, there is a balanced $r$-colouring of $K_n$ such that, in every copy of a $\lceil 1/\eps\rceil$-caterpillar on $n$ vertices in $K_n$, the number of edges in colour $1$ is smaller than the number of edges in one of the colours $2,\ldots,r$.
\end{proposition}

\noindent
While the above result only holds for a large number of colours, it turns out that a similar result holds for every $r\ge 3$ and \emph{all} trees with maximum degree at least $23n/24+O(1)$.

\begin{proposition}\label{prop:k>2}
	There is a constant $C>0$ such that the following holds.
 For every integer $r\ge 3$ and sufficiently large $n$, there is a balanced $r$-colouring of $K_n$ such that for every tree $T$ on $n$ vertices with maximum degree $\Delta \ge 23n/24+C$ and every copy of $T$ in $K_n$, the number of edges in colour $1$ is at most the number of edges in some colour among $2,\ldots,r$.
\end{proposition}

Propositions~\ref{prop:example_many_colours} and~\ref{prop:k>2} make it clear that additional conditions are necessary if we aim to find a copy of a tree $T$ in an arbitrary balanced $r$-colouring of $K_n$ where colour 1 is met substantially more often than any other colour.
A natural restriction, consistent with the previous results, is to impose an upper bound on the maximum degree of $T$.
However, even for spanning paths, certain almost balanced colourings can still behave pathologically.

\begin{example}
	\label{ex:1}
	Fix $r=4$, $n = 2N$ and let $V_1,V_2$ be a partition of the vertices of $K_n$ into two equal parts. 
	Colour all edges within $V_1$ (resp.\ $V_2$) in colour $1$ (resp.\ colour $2$), and assign colour $3$ (resp.\ colour $4$) to $\lfloor N^2/2\rfloor$ edges (resp.\ $\lceil N^2/2\rceil$ edges) between $V_1$ and $V_2$.
	In this almost balanced 4-colouring of $K_n$, every Hamilton cycle contains the same number of edges in colours $1$ and $2$.
	As a result, the difference between the number of edges in colours $1$ and $2$ in every Hamilton path is as most $1$.
\end{example}

Curiously, the construction in Example~\ref{ex:1} turns out to be the unique obstruction in a certain sense.
To make this statement formal, we introduce the notion of an $\eta$-standard edge-colouring of a graph $G$ on $n$ vertices.

\begin{definition}
	\label{def:StandardColouring}
	For every $\eta > 0$, an $r$-colouring of $G$ is \emph{$\eta$-standard} if, for every partition $V_1,V_2$ of its vertex set with $|V_1|-|V_2|\in \{0,1\}$ and every integer $i\in [2,r]$, $G[V_1]$ contains at least $\eta n^2$ edges with colour different from $1$ or $G[V_2]$ contains at least $\eta n^2$ edges with colour different from $i$.
\end{definition}
In essence, $\eta$-standard colourings are those that are ``not too close'' to the colouring described in Example~\ref{ex:1}, although the second colour is not necessarily colour 2. 
In particular, it is not hard to check that all balanced colourings in $r\neq 4$ colours are $\eta$-standard for all sufficiently small~$\eta$.

For any tree $T$ on $n$ vertices with suitably bounded maximum degree, the next theorem ensures that, in any almost balanced $\eta$-standard colouring of $K_n$, there is a copy of $T$ where colour 1 is met substantially more frequently than any other colour.
More precisely, we say that an $r$-colouring of a graph $G$ is \emph{$\eps$-balanced} if each colour $i$ is given to a $(1/r\pm \eps)$-proportion of the edges of $G$, that is, a proportion in the interval $[1/r - \eps, 1/r + \eps]$.

\begin{theorem}\label{thm:col_disc_2}
		Fix an integer $r\ge 3$ and a sufficiently small $\eta > 0$ (in terms of $r$). Then, there is $\eps = \eps(r,\eta) > 0$ such that, for every sufficiently large $n$ and every tree $T$ on $n$ vertices and maximum degree $\Delta\le \eps n$, every $\eps$-balanced $\eta$-standard colouring of $K_n$ in $r$ colours contains a copy of $T$ where colour $1$ is given to at least $(\eta^2/10^6 r^3) n$ more edges than any other colour.
\end{theorem}

It is natural to wonder if Theorem~\ref{thm:col_disc_2} extends to dense graphs in place of complete graphs.
The following example is inspired by a construction in~\cite{BCJP20} and shows that Example~\ref{ex:1} is not the only exceptional construction in case the minimum degree of the host graph drops below $3n/4$. 

\begin{example}\label{ex:2}
Let $n=4N$ and partition the vertex set of $K_n$ into two parts $V_1$ and $V_2$ of size $N$ and $3N$, respectively. 
Let $G$ be a graph obtained from deleting all edges connecting two vertices of $V_1$ as well as the edges of an $(N-1)$-regular graph on the vertex set $V_2$ (note that this is possible since $3(N-1)N$ is even). 
In particular, $G$ is a $3N$-regular graph.
Then, colour all edges incident to $V_1$ in blue and all other edges in red.
One can easily check that a spanning path of $G$ contains either $2N-2$ blue edges (if it starts and ends in $V_1$), or $2N-1$ blue edges (if it starts in $V_1$ and ends in $V_2$), or $2N$ blue edges (if it starts and ends in $V_2$), so the difference between the number of red and the number of blue edges is bounded independently of $N$ in each case.
\end{example}

While it would definitely be of interest to determine if $3n/4$ is a barrier for the existence of obstructions qualitatively different from Example~\ref{ex:1}, 
our approach allows us to treat host graphs of minimum degree at least $(1-\eps)n$ for small $\eps > 0$ depending on $r$, $\eta$ and the maximum degree of the tree.

\begin{theorem}\label{thm:col_disc_3}
		Fix integers $r\ge 3$, $\Delta\ge 2$ and a sufficiently small $\eta$ (in terms of $r$). Then, there is $\eps = \eps(r,\Delta,\eta) > 0$ such that, for every sufficiently large $n$, every graph $G$ on $n$ vertices and minimum degree at least $(1-\eps)n$, and every tree $T$ on $n$ vertices and maximum degree at most $\Delta$, every $\eps$-balanced $\eta$-standard colouring of $G$ in $r$ colours contains a copy of $T$ where colour $1$ is given to at least $\eps n$ more edges than any other colour.
\end{theorem}

We note that our proof techniques are sufficiently flexible to derive analogues of our main results for graphs containing many leaves or many adjacent pairs of vertices of degree two (such as sparse random graphs) instead of trees.
Nevertheless, we focus on trees for the sake of cleaner statements.

\subsection{\texorpdfstring{$k$-sum subgraphs of dense $2$-edge-coloured graphs}{k-sum subgraphs of dense 2-edge-coloured graphs}}
\label{subsubsec:k-sum-intro}
Given a suitably dense $2$-edge-coloured host graph, the rest of this work is dedicated to finding copies of a graph $H$ with a specified number of edges in each colour. 
For a $2$-edge-colouring $\bc:E(G)\to \{-1,1\}$ (where colour $-1$ is often called red, and colour $1$ is called blue) and a subgraph $G'\subseteq G$, we define $\bc(G') = \sum_{e\in E(G')} \bc(e)$.
A subgraph $G'\subseteq G$ is called \emph{$k$-sum} if $\bc(G') = k$.

For a $k$-sum copy of a graph $H$ in a host graph $G$ to exist, several elementary conditions must be satisfied: 
\begin{itemize}
    \item on the parity and the absolute value of $k$ ($k\equiv |E(H)|\, (\hspace{-0.7em}\mod 2)$ and $|k|\le |E(H)|$),
    \item on the structure of the graph $H$ and the host graph ($H$ needs to be realised as a subgraph of $G$ in the first place),
    \item and on the balancedness of $\bc$ (for example, if $|k|\neq |E(H)|$, no $k$-sum copies of $H$ exist in monochromatic host graphs).
\end{itemize}  
Taking these restrictions into account, it is natural to ask if, for a fixed $k$, there is always a $k$-sum copy of a graph $H$ with bounded maximum degree and a number of edges respecting the parity of $k$ in approximately balanced colourings of sufficiently dense host graphs.
We begin with a couple of motivating examples showing that this is not the case due to an additional number-theoretic constraint.

\begin{example}\label{ex:cycle_bip}
	Fix an integer $n = 4m+6$ and partition the vertices of a graph $K_n$ into two sets $V_1$ and $V_2$ such that $|V_1|=2m+1$ and $|V_2|=2m+5$. 
	Colour all edges between $V_1$ and $V_2$ in red and all other edges in blue. Then, every Hamilton cycle in $K_n$ has a discrepancy congruent to $2\,(\hspace{-0.7em}\mod 4)$.
    Indeed, given a Hamilton cycle of $K_n$, the number of red (resp.\ blue) edges in this colouring is $2\ell$ (resp.\ $n-2\ell$) for some integer $\ell\ge 0$, so its discrepancy is always congruent to $n$ modulo $4$.
\end{example}

Before stating the second example, we remark that it refutes the following conjecture of Mohr, Pardey and Rautenbach~\cite{mohr2022zero}.

\begin{conjecture}[Conjecture 1.5 in~\cite{mohr2022zero}]\label{conj:MPR}
	Fix integers $k,n\ge 1$ such that $\tbinom{n}{2}$ and $(k-1)n/k$ are both even integers. If $T$ is a tree of order $k$, $\bc : E(K_n) \to \{-1, 1\}$ is a 0-sum labeling of $K_n$, and $n$ is sufficiently large
	in terms of $k$, then $K_n$ has a 0-sum $T$-factor, that is, there is a 0-sum spanning forest $F$ of $K_n$
	whose components are all isomorphic to $T$.
\end{conjecture}

\begin{example}\label{ex:forest_bip}
	Fix $n=16k^2$ and $m = 8k^2-2k$ for some odd $k$. Partition the vertices of $K_n$ into sets $V_1$ and $V_2$ such that $|V_1|=m$ and $|V_2|=n-m$.
    Colour all edges between $V_1$ and $V_2$ in red and all other edges in blue. 
    In particular, the number of red edges is
	\[m(n-m) = (8k^2-2k)(8k^2+2k) = 64k^4 - 4k^2 = \frac{1}{2}\binom{n}{2}.\]
	Let $S$ be a star on $8$ vertices. We show that no $S$-factor in $K_n$ is 0-sum.
	Fix an $S$-factor $H$ and let $H_1$ and $H_2$ be the subforests containing all stars with centers in $V_1$ and $V_2$, respectively. 
	Let $s_1,s_2$ be the number of stars in $H_1,H_2$ and $\ell_1,\ell_2$ be the number of red edges in $H_1,H_2$, respectively.
	Then, $8s_1-\ell_1+\ell_2 = m$, so $\ell_1+\ell_2$ is even. Moreover, $s_1+s_2=n/8=2k^2$ and, since $k$ is odd, we obtain that $s_1+s_2 \equiv 2(\hspace{-0.7em} \mod 4)$.
	As a result,
	\[\bc(H) = (m-s_1-\ell_2) + (n-m-s_2-\ell_1) - (\ell_1+\ell_2) = n - 2(\ell_1+\ell_2) - (s_1+s_2) \equiv 2(\hspace{-1em} \mod 4).\]
 This shows that there is no 0-sum $S$-factor in this graph.
\end{example} 

In fact, the absence of 0-sum copies in these examples reflects a more general phenomenon. Its essence is captured by the following definition and the subsequent proposition.

\begin{definition}\label{def:pair}
	For an integer $n\ge 1$ and a spanning subgraph $G$ of $K_n$, a colouring $\bc: E(K_n)\to \{-1,1\}$ is called \emph{bipartite} if there is a partition $(V_1,V_2)$ of the vertices of $K_n$ such that either $\bc^{-1}(1)$ or $\bc^{-1}(-1)$ is equal to the edge set of $G[V_1]\cup G[V_2]$.
	Then, the sets $V_1,V_2$ are called the \emph{parts} of $\bc$.
	Moreover, for an $n$-vertex graph $H$ and $G, \bc$ as above, the triplet $(H,G,\bc)$
	is called \emph{inert} if all degrees in $H$ have the same parity and the colouring $\bc$ is bipartite.
\end{definition}

Note that, in each of Example~\ref{ex:cycle_bip} and Example~\ref{ex:forest_bip}, $(H, K_n, \bc)$ is an inert triplet.
This is not a coincidence, as all inert triplets share the following property.

\begin{proposition}\label{prop:par-bip}
	For every $n\ge 1$, every inert triplet $(H, G, \bc)$ and every copy $H'$ of $H$ in $G$, $\bc(H')$ has the same residue modulo $4$.
\end{proposition}

More interestingly, inertness turns out to be the only obstacle for finding $k$-sum copies of a graph subject to some natural conditions, as demonstrated by the following theorem. 

\begin{theorem}
	\label{thm:AlmostZeroSumDenseGraph}
	For every integer $\Delta\ge 1$, there are constants $\eps, \nu, \alpha > 0$ such that the following holds for all sufficiently large $n$ and $m$.
	Let $H$ be an $n$-vertex graph with $m$ edges and maximum degree $\Delta$, let $G$ be an $n$-vertex graph with minimum degree at least $(1-\eps)n$, and let $\bc$ be a $\nu$-balanced $2$-edge-colouring of $G$. 
	If the triplet $(H,G,\bc)$ is inert, let $a=4$; otherwise, let $a=2$.
	Then, there is an integer $b\in \{0,1,2,3\}$ such that, for every integer $k\in [-\alpha m, \alpha m]$, there is a copy $H'\subseteq G$ of $H$ with $\bc(H') = k$ if and only if $k\equiv b\,(\hspace{-0.62em}\mod a)$.
\end{theorem}

Note that, in the case of a non-inert triplet $(H,G,\bc)$, \Cref{thm:AlmostZeroSumDenseGraph} ensures the existence of a copy $H$ where the number of edges in the two colours differ by at most $2$. 
Before outlining our proofs, we note that the problem of finding spanning copies that are nearly colour-balanced has been addressed in several papers. The specific case of perfect matchings has been studied in \cite{ehard2020low,Hol24, kittipassorn2023existence,pardey2022almost}, while factors of $P_3$'s and $P_4$'s were considered in~\cite{mohr2022zero}, and spanning paths were the object of~\cite{caro2022zero}. 
The case of spanning forests of $K_n$ with fixed isomorphism class was investigated in~\cite{hollom2024almost,pardey2023efficiently}.

\paragraph{Brief outline of the proofs.} The proofs of Theorems~\ref{thm:col_disc_1},~\ref{thm:col_disc_2} and~\ref{thm:col_disc_3} rely on a switching technique that allows us to increase or decrease the number of edges in a copy of $T$ with a given colour via small local changes.
More precisely, given a tree $T$ on $n$ vertices randomly embedded in the complete graph $K_n$, we find some special pairs of edges $uv, wz$ in $T$ (described in the beginning of Section~\ref{sec:thm1}) which can be exchanged with the edges $uz, vw$ outside $T$ such that the resulting graph remains a copy of $T$.
In the proof of Theorem~\ref{thm:col_disc_1}, these switchings are used solely to increase the number of edges of a fixed colour.
However, a more careful treatment is needed in the proof of Theorem~\ref{thm:col_disc_2} where we need to show that, for every $i\in [r]\setminus \{1\}$, a significant number of switchings increase the number of edges in colour 1 but not the number of edges in colour $i$.
As pointed out in Example~\ref{ex:1}, this is not always possible; however, by counting suitably coloured $4$-cycles in edge-colourings of $K_n$, we show that this example is the only obstruction.
Additional complications arise in Theorem~\ref{thm:col_disc_3} where edges outside the host graph cannot be used.
This proof can be seen as a modified version of the argument showing Theorem~\ref{thm:col_disc_2} where some additional vertex exchanges are used to avoid the non-edges of $G$.
The construction showing Proposition~\ref{prop:example_many_colours} uses a colouring of a blow-up of a small complete graph and is inspired by the finite projective plane $\mathrm{PG}(2,p)$.
In contrast, Proposition~\ref{prop:k>2} employs a more straightforward colouring of the same blow-up.

We turn to the complementary $k$-sum part.
The proof of \Cref{thm:AlmostZeroSumDenseGraph} combines ideas from previous sections with new structural insights.
First, unless $(H,G,\bc)$ is an inert triplet, we find and reserve a small part of the graph $G$ (later called a $2$-gadget), which will enable us to adjust the sum of the graph by exactly $+2$ or $-2$.
We then embed the graph $H$ uniformly at random in $K_n$. 
We show that, with positive probability, this random embedding contains \begin{itemize}
    \item $\Omega(m)$ disjoint sets of vertices that allow us to change $\bc(H)$ by $\pm 4$ (called gadgets),
    \item a small linear number of edges outside $G$,
    \item and close to the expected number of edges in each colour.
\end{itemize}
Then, we exchange some pairs of vertices to obtain some specific configuration in the reserved part of $G$ and to remove any missing edge in the embedding. 
This process is carefully executed to ensure that
\begin{itemize}
	\item most of the gadgets remain untouched, and
	\item the number of edges in each colour is not altered by much.
\end{itemize}
We finish the proof of \Cref{thm:AlmostZeroSumDenseGraph} by triggering some of the untouched gadgets to produce an embedding of sum precisely $k$. 
Finally, the proof of Proposition~\ref{prop:par-bip} relies on a simple degree-counting argument in and between the parts of a bipartite colouring. 

\paragraph{Plan of the paper.} 
The remainder of this paper is organised as follows.
In \Cref{sec:prem}, we describe the switching method and prove several useful lemmas.
In \Cref{sec:thm1}, we prove Theorem~\ref{thm:col_disc_1} and deduce \Cref{col:discErdos}.
\Cref{sec:3} is dedicated to the proofs of \Cref{prop:example_many_colours} and~\ref{prop:k>2}, and \Cref{sec:5} deals with the proofs of \Cref{thm:col_disc_2} and \Cref{thm:col_disc_3}. 
Then, in Section~\ref{sec:5'}, we prove Proposition~\ref{prop:par-bip} and Theorem~\ref{thm:AlmostZeroSumDenseGraph}.
We conclude the paper with additional remarks and open questions in Section~\ref{sec:conc}.

\section{Preliminary results and tools}
\label{sec:prem}

In the proofs in this and the following sections, upper and lower integer parts will be omitted when irrelevant for the argument.

\subsection{Switchable edges and the switching technique}
Fix a tree $T$ on $n$ vertices. A pair of vertex-disjoint edges $uw$ and $vz$ is called \emph{switchable} if one of the following conditions holds.
\begin{itemize}
	\item Both edges contain a leaf of $T$, say $w$ in $uw$ and $z$ in $vz$ (switchable pair of type I).
	\item Both edges contain a vertex of degree two, say $u$ in $uw$ and $v$ in $vz$. Moreover, $u$ and $v$ share a common neighbour different from $w$ and $z$ (switchable pair of type II).
	\item Both edges contain a vertex of degree two, say $u$ in $uw$ and $v$ in $vz$. Moreover, $uv$ is an edge in $T$ (switchable pair of type III).
\end{itemize}

\begin{figure}[ht]
	\centering
	\begin{tikzpicture}
    % Common settings for nodes to reuse
    \tikzset{vertex/.style={circle, fill, inner sep=1.5pt}}
    
    % First diagram
    \begin{scope}[xshift=0cm]
        % Vertices
        \node[vertex, label=left:$w$] (u1) at (0, 0) {};
        \node[vertex, label=right:$z$] (v1) at (2, 0) {};
        \node[vertex, label=left:$u$] (w1) at (0, 2) {};
        \node[vertex, label=right:$v$] (z1) at (2, 2) {};
        
        % Edges
        \draw (u1) -- (w1);
        \draw (v1) -- (z1);

        % Circles
        \draw[thick] (u1) circle (0.15);
        \draw[thick] (v1) circle (0.15);
        
        % Label
        \node at (1, -0.75) {I};
    \end{scope}
    
    % Second diagram
    \begin{scope}[xshift=5cm]
        % Vertices
        \node[vertex, label=left:$w$] (u2) at (0, 0) {};
        \node[vertex, label=right:$z$] (v2) at (2, 0) {};
        \node[vertex, label=left:$u$] (w2) at (0, 2) {};
        \node[vertex, label=right:$v$] (z2) at (2, 2) {};
        \node[vertex] (y2) at (1,3) {};
        
        % Edges
        \draw (u2) -- (w2);
        \draw (v2) -- (z2);
        \draw (w2) -- (y2) -- (z2);

        % Circles
        \draw[thick] (w2) circle (0.15);
        \draw[thick] (z2) circle (0.15);
        
        % Label
        \node at (1, -0.75) {II};
    \end{scope}
    
    % Third diagram
    \begin{scope}[xshift=10cm]
        % Vertices
        \node[vertex, label=left:$w$] (u3) at (0, 0) {};
        \node[vertex, label=right:$z$] (v3) at (2, 0) {};
        \node[vertex, label=left:$u$] (w3) at (0, 2) {};
        \node[vertex, label=right:$v$] (z3) at (2, 2) {};
        
        % Edges
        \draw (u3) -- (w3);
        \draw (v3) -- (z3);
        \draw (w3) -- (z3);

        % Circles
        \draw[thick] (w3) circle (0.15);
        \draw[thick] (z3) circle (0.15);
        
        % Label
        \node at (1, -0.75) {III};
    \end{scope}
\end{tikzpicture}
    \caption{An illustration of a switchable pair of edges $uw, vz$ of type I, II and III, respectively. The vertices marked with small circles have all their neighbours in the figure.}
	\label{fig:switchable-edges}
\end{figure}
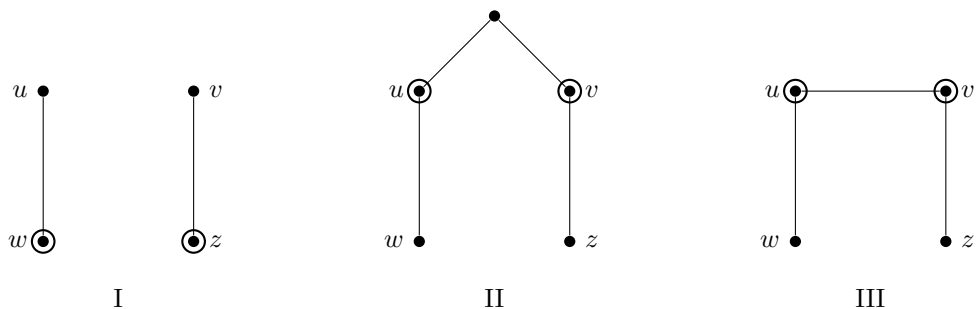

Note that, for each of the three types of switchable pairs, the switching $uw,vz\to uz,vw$ leaves us with a tree isomorphic to $T$. We call this switching \emph{allowed}, see Figure~\ref{fig:switchable-edges}.
Observe that a switching is an \emph{operation} of replacing a pair of edges in some embedding of $T$ with another pair of edges, thus forming another embedding of the tree, while the term switchable is a \emph{characteristic} of a pair of edges indicating that the switching operation applied to them results into a tree isomorphic to $T$.

We say that a family $\mathcal S$ of switchable edge pairs is \emph{free} if all pairs in it are edge-disjoint.
In particular, this means that every allowed switching of a pair in $\mathcal S$ can be done in parallel of all other allowed switchings of pairs in $\mathcal{S}$.
We also say that a path $P$ in $T$ is a \emph{bare path} if every vertex in $P$ except for possibly the endpoints has degree two in $T$.

\begin{lemma}\label{lem:many switches}
	Let $T$ be a tree on $n$ vertices with maximum degree $\Delta$.
	Then $T$ contains a free family of pairs of switchable edges of size at least $(n-1-\Delta)/40$.
\end{lemma}
\begin{proof}
	Set $\zeta = \zeta(n,\Delta) = \lceil (n-1-\Delta)/40\rceil$.
	First, we consider three special cases. 
	\begin{enumerate}
		\item If $T$ is a path of length $n-1\ge 3$, one can find $\lfloor (n-1)/3\rfloor\ge \zeta$ edge-disjoint paths of length three, each containing a switchable pair of type III.
		Moreover, these switchable pairs are pairwise disjoint, so we get a free family of the required size. 
		\item If $\zeta = 0$, the statement is trivially satisfied.
		\item If $\zeta = 1$, then $T$ contains at least one pair of leaves with different parents that form an allowed switching of type I.
	\end{enumerate}
	
	In the remainder of the proof, we assume that $\Delta\in [3,n-42]$ (so $n\ge 45$).
	For every $i\in [\Delta]$, denote by $d_i$ the number of vertices of degree precisely $i$ in $T$.
	By double counting the vertices and the edges of $T$, respectively, we get that
	\begin{equation}\label{eq:tree_relations}
		\sum_{i=1}^{\Delta} d_i = n\quad \text{and}\quad \sum_{i=1}^{\Delta} id_i = 2(n-1).
	\end{equation}
	We consider several cases.
	
	\paragraph{Case 1: $d_1\ge 2\zeta - 1 + \Delta$.} Let $\mathcal{L}$ be the family of edges of $T$ containing a leaf.
	Then, $|\mathcal{L}| = d_1$ and for every edge $e$ in $\mathcal{L}$, there are at most $\Delta-1$ other edges in $\mathcal{L}$ that do not form a switchable pair of type I with $e$.
	Hence, by forming disjoint switchable pairs greedily, we end up with at least $\lceil (d_1 - \Delta)/2\rceil\ge \zeta$ such pairs, as desired.
	
	\paragraph{Case 2: $d_1 \le 2\zeta - 2 + \Delta$.} Note that~\eqref{eq:tree_relations} implies $\sum_{i=1}^{\Delta} 2d_i = 2+\sum_{i=1}^{\Delta} id_i$.
	Recalling that $\Delta \ge 3$, we obtain
	\begin{align*}
		d_1 &= 2+\sum_{i=3}^{\Delta} (i-2)d_i = \Delta + \bigg(\sum_{i=3}^{\Delta} (i-2)d_i\bigg) - (\Delta-2) \\
		&= \Delta + \sum_{i=3}^{\Delta}d_i + \sum_{i=4}^{\Delta-1}(i-3)d_i + (\Delta-3)(d_{\Delta}-1) - 1 \ge \Delta + \sum_{i=3}^{\Delta}d_i - 1
	\end{align*}
	By the assumption on $d_1$ this yields that the number of vertices of degree at least three satisfies
    \begin{align}
    \label{eq:SumdiUB}
        \sum_{i=3}^{\Delta} d_i \le d_1 + 1 - \Delta \le 2\zeta -1,
    \end{align}
	implying that the number of vertices of degree two in $T$ satisfies
     \begin{align}
     \label{eq:d2LB}
         d_2 \ge n - (2\zeta - 2 +\Delta) - (2\zeta - 1) = n - \Delta - 4\zeta + 3.
     \end{align}
	
	Now, for every integer $i\ge 1$, denote by $\ell_i$ the number of bare paths in $T$ consisting of exactly $i$ internal vertices of degree two, and with its endpoints being of degree different than two.
	In particular, $\sum_{i=1}^{\infty} i\ell_i = d_2$.
	We further consider two subcases.
	
	\paragraph{Case 2.1: $\sum_{i=1}^{\infty} \ell_i \ge 4\zeta-2$.} 
    Note that, given two bare paths containing the subpaths $uw_1z_1$ and $uw_2z_2$ (with $w_1, w_2$ of degree two in $T$), the edges $w_1z_1$ and $w_2z_2$ form a switchable pair of type II. 
	Thus, we greedily form a free family of switchable pairs of type II by consecutively pairing bare paths that share a common endvertex. 
	Note that every bare path counted by the $\ell_i$'s must contain at least one vertex of degree greater than two (recall that $T$ is not a bare path itself).
    Therefore, as long as there are at least $2\zeta$ yet unpaired bare paths, at least two of them must share an endvertex of degree greater than two, as there are at most $2\zeta - 1$ vertices of degree at least three by \eqref{eq:SumdiUB}.
	Hence, these two paths can be paired, and this process ensures that at least $\lceil \left(\sum_{i=1}^{\infty}\ell_i - (2\zeta-1) \right)/2\rceil\ge \zeta$ pairs of bare paths can be formed, as desired.
	
	\paragraph{Case 2.2: $\sum_{i=1}^{\infty} \ell_i\le 4\zeta -3$.} Then, given a bare path $P$ containing at least $i$ vertices of degree two for some $i\ge 1$, one can find $\lfloor\frac{i+1}{3} \rfloor$ edge-disjoint subpaths of length three in $P$.
	Moreover, the first and the last edge in each of these subpaths form a switchable pair of edges of type III. 
	Thus, the number of disjoint switchable edge pairs is bounded from below by
	\[\sum_{i=1}^{\infty} \left\lfloor\frac{i+1}{3} \right\rfloor\ell_i\ge \sum_{i=1}^{\infty} \frac{i-1}{3} \ell_i = \frac{1}{3}\left(d_2 - \sum_{i=1}^{\infty}\ell_i \right) \ge \frac{1}{3}\left(n - \Delta - 8\zeta + 6 \right) \ge \zeta,\]
	where the second inequality follows from \eqref{eq:d2LB} and the last inequality from the definition of $\zeta$ and $n\ge \Delta+42$.
	Since the above switchable pairs of type III form a free family, the proof is completed.
\end{proof}

\subsection{Marked \texorpdfstring{$4$}{4}-cycles}\label{sec:gen_marked}

We define a \emph{marked $4$-cycle} to be a cycle $\gamma\subseteq K_n$ of length 4 where one of the two pairs of opposite edges in the cycle is said to be \emph{marked}.
These marked edges will later serve us to indicate which of the two possible switchings along the cycle is advantageous for us. 
However, we keep the setting general to allow for more flexibility and potentially further applications of this approach.

Consider a switchable pair of edges $uw,vz$ in $T$, an embedding of $T$ in $K_n$ and a family $\mathcal{C}$ of marked $4$-cycles in $K_n$.
The switching $uw, vz\to uz, vw$ is called \emph{desirable} (with respect to $\mathcal{C}$) if these edges are embedded onto a pair of marked edges in a cycle $\gamma\in \mathcal{C}$ so that (the images of) $u$ and $v$ are at distance two in $\gamma$.
The next lemma shows that, given families $\mathcal{C}_1, \mathcal{C}_2, \ldots, \mathcal{C}_k$ of marked 4-cycles, each of size at least $c n^4$, a uniformly random embedding of a tree $T$ with sufficiently small maximum degree in $K_n$ typically contains linearly many desirable switchings with respect to each family.

\begin{lemma}\label{lem:marked_4_cycles}
	Fix constants $\eps\in (0,1/2)$ and $c\in (0,1/4)$, an integer $k\ge 1$ and families $\mathcal{C}_1, \mathcal{C}_2, \ldots, \mathcal{C}_k$, each containing at least $c n^4$ marked $4$-cycles in $K_n$. 
	Then, for every sufficiently large $n$ and every tree $T$ on $n$ vertices with maximum degree $\Delta\le \eps n$, the following holds with probability at least $1-400\eps k/c$.
	In a uniformly random embedding of $T$ in $K_n$, there exists a free family containing at least $cn/40$ switchable pairs in desirable switchings with respect to $\mathcal{C}_i$ for every~$i\in [k]$.
\end{lemma}
\begin{proof}
	Fix a tree $T$ with maximum degree $\Delta\le \eps n$ and identify it with its random embedding in $K_n$.
	Then, Lemma~\ref{lem:many switches} ensures the existence of a free family $\mathcal S$ containing $s \defined |\mathcal S| \ge (n-1-\Delta)/40 \ge n/80$ switchable pairs.
	We show via a second moment computation that, for every $i\in [k]$, up to choosing~$\eps$ small, typically a constant proportion of the switchable pairs in $\mathcal S$ form a desirable switching with respect to $\mathcal C_i$.
	
	Given $i \in [k]$, denote by $\alpha_i$ the probability that a marked $4$-cycle chosen uniformly at random in $K_n$ is in~$\mathcal C_i$.
	As there are $n(n-1)(n-2)(n-3)/4$ possible marked $4$-cycles in $K_n$, we have
	\[\alpha_i \coloneqq \frac{4|\mathcal C_i|}{n(n-1)(n-2)(n-3)} \ge 4c. \]
	Let $X_i$ be the random variable counting the number of pairs of edges in $\mathcal S$ such that their embeddings in $K_n$ form a desirable switching with respect to $\mathcal C_i$.
	Given a pair of edges $(e,f) \in \mathcal S$ denote by $I_{(e,f)}$ the indicator random variable of the event that $(e,f)$ participate in a desirable switching with respect to $\mathcal C_i$.
	Then,
	\[X_i = \sum_{(e,f) \in \mathcal S} I_{(e,f)}\quad \text{and}\quad \mathbb E\left[X_i \right] = s \alpha_i \ge 4c\cdot \frac{n}{80} = \frac{cn}{20}.\]

We turn to a computation of $\mathrm{Var}(X_i)$. 
Note that, for every pair of edges $(e,f)\in \cS$, there are at most $4\Delta\le 4\eps n$ pairs of edges $(e',f')\in \cS$ such that $e\cup f$ and $e'\cup f'$ share a common vertex. 
Moreover, for all vertex-disjoint pairs of edges $(e,f), (e',f') \in \mathcal S$, the probability that both $(e,f)$ and $(e',f')$ are embedded into a desirable switching with respect to $\mathcal C_i$ is at most
	\[\frac{16|\mathcal C_i|^2}{\prod_{i=0}^{7} (n-i)} = \alpha_i^2 \cdot \frac{\prod_{i=0}^{3} (n-i)}{\prod_{i=4}^{7} (n-i)} = \alpha_i^2 \cdot \frac{\prod_{i=4}^{7} (n-i+4)}{\prod_{i=4}^{7} (n-i)} \le \alpha_i^2\left(1 + \frac{20}{n} \right). \]
	Hence, by using that $s\alpha_i\ge cn/20$, we obtain
	\begin{align}
	\begin{split}
	\label{eq:var}
		\Var\left(X_i \right) &=  \sum_{(e,f),(e',f') \in \mathcal S} \left(\mathbb E\left[I_{(e,f)} I_{(e',f')} \right] - \mathbb E\left[I_{(e,f)} \right] \mathbb E\left[I_{(e',f')} \right] \right) \\
        &\le \sum_{(e,f)\in \mathcal S} 4\eps n \mathbb E\left[I_{(e,f)}\right] + s^2 \alpha_i^2 \cdot \frac{20}{n} =  4\eps n\cdot s\alpha_i + s^2 \alpha_i^2 \cdot \frac{20}{n}.
	\end{split}
	\end{align}

    Finally, by Chebyshev's inequality,
	\[\mathbb P(X_i \le cn/40) \le \mathbb P(X_i\le \mathbb E[X_i]/2)\le \frac{4\mathrm{Var}(X_i)}{\mathbb E[X_i]^2} \le \frac{16\eps n\cdot s\alpha_i}{s^2\alpha_i^2} + \frac{20}{n}\le \frac{400\eps}{c},\]
	and a union bound over all $i\in [k]$ finishes the proof.
\end{proof}

\subsection{\texorpdfstring{$i$-marked $4$}{4}-cycles in edge-coloured graphs}\label{sec:4.2}

\begin{figure}[ht]
	\centering
	\begin{tikzpicture}
    % Common settings for nodes to reuse
    \tikzset{vertex/.style={circle, fill, inner sep=1.5pt}}
    
    % First diagram
    \begin{scope}[xshift=0cm]
        % Vertices
        \node[vertex, label=left:$x$] (x1) at (0, 0) {};
        \node[vertex, label=left:$y$] (y1) at (0, 2) {};
        \node[vertex, label=right:$u$] (u1) at (2, 2) {};
        \node[vertex, label=right:$v$] (v1) at (2, 0) {};
        
        % Edges
        \draw[blue] (x1) -- (y1) -- (u1) -- (v1);
        \draw[dashed, black] (v1) -- (x1);
        
        % Label
        \node at (1, -0.75) {(i)};
    \end{scope}
    
    % Second diagram
    \begin{scope}[xshift=5cm]
        % Vertices
        \node[vertex, label=left:$x$] (x2) at (0, 0) {};
        \node[vertex, label=left:$y$] (y2) at (0, 2) {};
        \node[vertex, label=right:$u$] (u2) at (2, 2) {};
        \node[vertex, label=right:$v$] (v2) at (2, 0) {};
        
        % Edges
        \draw[blue] (x2) -- (y2);
        \draw[blue] (u2) -- (v2);
        \draw[dashed, black] (v2) -- (x2);
        \draw[dashed, black] (y2) -- (u2);
        
        % Label
        \node at (1, -0.75) {(ii)};
    \end{scope}
    
    % Third diagram
    \begin{scope}[xshift=10cm]
        % Vertices
        \node[vertex, label=left:$x$] (x3) at (0, 0) {};
        \node[vertex, label=left:$y$] (y3) at (0, 2) {};
        \node[vertex, label=right:$u$] (u3) at (2, 2) {};
        \node[vertex, label=right:$v$] (v3) at (2, 0) {};
        
        % Edges
        \draw[blue] (x3) -- (y3);
        \draw[dashed, black] (y3) -- (u3);
        \draw[dashed, black] (v3) -- (x3);
        \draw[dashed] (u3) --  (v3);
        
        % Label
        \node at (1, -0.75) {(iii)};
    \end{scope}
\end{tikzpicture}
	\caption{The three cases of a marked 4-cycle. Solid edges are blue, and dashed edges are some colour other than blue.} 
	\label{fig:marked-cycles}
\end{figure}
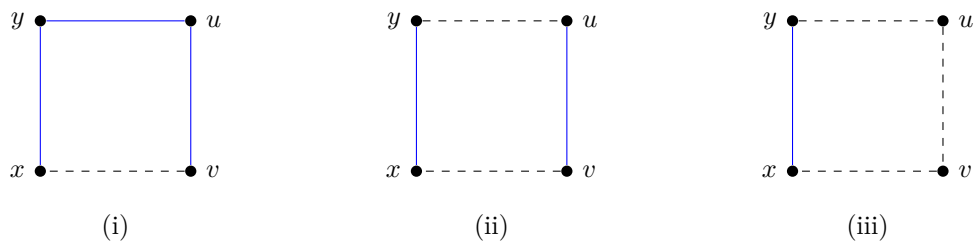

In this section, we restrict the general setting introduced in Section~\ref{sec:gen_marked}.
Fix an $r$-edge-colouring $E(K_n)\mapsto [r]$.
For our purposes, a 4-cycle will be marked if the corresponding desirable switching increases the number of edges in colour 1, which we shall refer to as \emph{blue} for convenience.
More precisely, a 4-cycle $\gamma = xyuv$ is \emph{marked} if one of the following holds (see~\Cref{fig:marked-cycles}):
\begin{enumerate}[(i)]
	\item $xy, yu, uv$ are blue and $vx$ is not blue, or
	\item $xy,uv$ are blue and $yu,vx$ are not blue, or
	\item $xy$ is blue and the other edges are not blue.
\end{enumerate}

Note that, in each of the three cases, the marked edges in $\gamma$ are $yu$ and $vx$.
Moreover, in the last case, we call the cycle \emph{$i$-increasing} if $\gamma$ contains a single edge in colour $i\in [2,r]$ and this edge is opposite to the blue one.

With an eye on the proof of \Cref{thm:col_disc_2,thm:col_disc_3}, for every $i\in [2,r]$, we consider the family $\mathcal C_i$ consisting of the marked 4-cycles of types (i) and (ii) as well as all non-$i$-increasing 4-cycles.
The goal of this section is to show that, for a tree $T$ as in Theorem~\ref{thm:col_disc_2}, which is embedded uniformly at random in a nearly balanced edge-colouring of $K_n$, we can make a linear number of desirable switchings in such a way that there is no $i$ for which almost all of these desirable switchings are $i$-increasing.
In general this is not always possible, as~\Cref{ex:1} shows; however, the next lemma proves that such problematic scenarios do not happen under the assumption that the colouring is $\eta$-standard for some small $\eta > 0$.

\begin{lemma}\label{lem:not-red-increasing}
	Fix an integer $r\ge 2$, $c \le 2^{-39}r^{-5}$ and $\eps\in (0,c)$.
	Consider an $\eps$-balanced, $32r\sqrt{rc}$-standard $r$-colouring of the complete graph $K_n$, and let $G$ be a spanning subgraph of $K_n$ with minimum degree $\delta(G)\geq (1-\eps)n$.
	Then, for sufficiently large $n$ and for every colour $i\in [2,r]$, there are at least $c n^4$ marked $4$-cycles in $G$ that are not $i$-increasing.
\end{lemma}

\begin{proof}
	Let $\alpha \defined \sqrt{8rc}$ and fix a colour $i\in [2,r]$, denoting the edges in this as \emph{red}.
	Recalling that we refer to colour 1 as blue, a 4-clique in $G$ is called \emph{blue} if all 6 of its edges are blue, and \emph{good} if it contains two vertex-disjoint blue edges but is not a blue clique.
	Note that every good 4-clique contains a marked 4-cycle which is not $i$-increasing. 
	Our proof proceeds in several steps.\\
	
	\noindent
	\textbf{Step 1: an approximate blue clique.} Let $E_B$ be the set of blue edges $e$ for which there are at most $\alpha n^2$ blue edges $f$ disjoint from $e$ such that $e\cup f$ does not induce a blue clique in $G$. 
	Note that for every pair of edges $e,f$ as above, if the clique they span is included in $G$, then it is good. We consider two cases.
	
	If $|E_B|$ contains less than a $(1-\alpha)$-proportion of all blue edges, then at least $\alpha (1/r-\eps)\tbinom{n}{2}$ blue edges are outside $E_B$. 
	Moreover, there are at most $n\cdot \eps n\cdot \binom{n}{2}\le \eps n^4/2$ quadruplets of vertices that span a graph with fewer than six edges in $G$.
	Since every good clique is spanned by at most two pairs of disjoint blue edges, the number of good cliques (and, as a consequence, of non-red-increasing marked 4-cycles) in $G$ is at least
	\[\frac{1}{2}\cdot \left(\alpha\left(\frac{1}{r}-\eps\right)\binom{n}{2}\cdot \alpha n^2 - \frac{\eps}{2}n^4 \right)\ge c n^4,\]
	which finishes the proof in this case.
	
	From now on, we assume that at least a $(1-\alpha)$-proportion of the blue edges are in $E_B$. Let $B\sseq V(G)$ be the set of endpoints of blue edges in $E_B$.
	We will show that almost all edges induced by $B$ in $G$ are blue and almost all blue edges go between two vertices in $B$.
	Indeed, for every blue edge $e$ in $E_B$, there are at least
	\begin{align*}
		s = s(\eps,r) \defined \left(\frac{1}{r}-\eps\right) \binom{n}{2} - (2n-3) - 2\eps n\cdot n - \alpha n^2\geq \left(\frac{1}{2r} - 2\alpha\right) n^2
	\end{align*} 
	blue edges $f$ in $E_B$ disjoint from $e$ and such that $e\cup f$ induces a blue 4-clique in $G$. Here, the term $2n-3$ is an upper bound for the number of blue edges intersecting $e$ and $2\eps n\cdot n$ dominates the number of blue edges $f$ such that the graph spanned by $e\cup f$ is not a 4-clique in $G$.
	Then, for every edge $e\in E_B$, the edges $f\in E_B$ forming blue cliques with 
	$e$ span a graph with at least $\sqrt{2s}\geq (1-\zeta) n/\sqrt{r}$ vertices, where $\zeta \defined 1-(1-4r\alpha)^{1/2}$.
	Hence, every vertex in $B$ is incident to at least $(1-\zeta) n/\sqrt{r}$ blue edges leading towards other vertices in~$B$. 
	Moreover, there are at most $(1/r+\eps)n^2/2$ blue edges, so $|B|\cdot (1-\zeta)n/\sqrt{r}\le (1/r+\eps)n^2$. Thus,
	\begin{equation}\label{eq:bounds_B}
		(1-\zeta) \frac{n}{\sqrt{r}}\le |B|\le (1+2\zeta) \frac{n}{\sqrt{r}},
	\end{equation}
	where for the second inequality we used the fact that $1+r\eps\le (1+2\zeta)(1-\zeta)$. Simply put,~\eqref{eq:bounds_B} shows that $E_B$ spans a slightly perturbed clique of size roughly $n/\sqrt{r}$ from~$G$. \\
	
	\noindent
	\textbf{Step 2: the case $r\in \{2,3\}$.} In this case, using that $0.55 < (1-\zeta)/\sqrt{3} < (1+2\zeta)/\sqrt{2} < 0.71$ together with~\eqref{eq:bounds_B} shows that $|B|\in [0.55n, 0.71n]$. Denote by $\mathcal B$ the family of 4-cycles in $K_n$ with two consecutive vertices in $B$ and two consecutive vertices outside $B$.
	Using that $x\in [1/2,1]\mapsto x(1-x)$ is a decreasing function, we get that
	\[|\mathcal B| = |B|(|B|-1)(n-|B|)(n-|B|-1)/2\ge 0.71n(0.71n-1)\cdot 0.29n(0.29n-1)/2 > n^4/100.\]
	Moreover, $\mathcal B$ contains at most $n\cdot \eps n\cdot n^2 = \eps n^4$ cycles containing a non-edge of $G$, at most $|B|\cdot (3\zeta n/\sqrt{r})\cdot n^2\le 3\zeta n^4$ cycles intersecting $B$ in a non-blue edge, and at most $\alpha(1/r+\eps)\tbinom{n}{2}\cdot n^2\le \alpha n^4$ cycles containing a blue edge with at most one endpoint in $B$. 
	Thus, when $r=2$, there must be at least $(1/100-\eps-3\zeta-\alpha)n^4 \ge cn^4$ cycles in $G$ with one blue and three red edges, which are in particular non-red-increasing, as desired.

    When $r=3$, call the third colour green. Since $|B|\cdot (n-|B|) > n^2/5\ge (1/3+\eps)\tbinom{n}{2} + n^2/50$, there must be at least $(n^2/50) (n^2/50-2n)\ge n^4/10^4$ cycles in $\mathcal B$ containing a non-green edge between $B$ and its complement. Applying similar reasoning to that used in the case $r=2$, we obtain at least $(1/10^4-\eps - 3\zeta-\alpha)n^4 \geq cn^4$ non-red-increasing cycles in $G$, which concludes the proof for $r=3$. \\ 
	
	\noindent
	\textbf{Step 3: the case $r\ge 4$.} In this case, our assumption that the colouring is $32r\sqrt{rc}$-standard comes into play. 
	Let $C$ be the complement of $B$.
	
	\begin{claim}
		\label{claim:many-not-red}
		For every $r\ge 4$, there are at least $16r\sqrt{rc}\cdot n^2$ non-red edges in $C$.
	\end{claim}
	\begin{proof}[Proof of Claim~\ref{claim:many-not-red}]
		If $r\geq 5$, the statement is immediate since there are at most $(1/5+\eps)\tbinom{n}{2}$ red edges and at least $\tbinom{n-|B|}{2} - \eps n^2\ge (1/5+\eps)\tbinom{n}{2}+16r\sqrt{rc}\cdot n^2$ edges of $G$ with two endpoints in $C$.
		
		Suppose that $r=4$. Let $B'$ be a set of size $\lfloor n/2\rfloor$ that either contains or is contained in $B$, and let $C'$ be the complement of $B'$. Given that our colouring is $32r\sqrt{rc}$-standard, there are either at least $32r\sqrt{rc}\cdot n^2$ non-blue edges in $B'$ or at least $32r\sqrt{rc}\cdot n^2$ non-red edges in $C'$.
		We will prove that the former option does not hold.
		On the one hand, using~\eqref{eq:bounds_B} we have that $|B'\setminus B|\le n/2 - (1-\zeta)n/2 = \zeta n/2$, and therefore, there are at most $\zeta n^2/4$ pairs of vertices in $B'$ with at least one vertex in $B'\setminus B$. 
		Moreover, using once again~\eqref{eq:bounds_B} and the fact that every vertex in $B$ is connected by blue edges to at least $(1-\zeta)n/2$ other vertices in $B$, we get that $B$ contains at most $|B|\cdot 3\zeta n/2\le 3\zeta n^2/2$ pairs of vertices that are not connected by a blue edge. 
		Summing the two bounds shows that there are at most $2\zeta n^2  < 32r\sqrt{rc}\, n^2$ pairs of vertices in $B'$ not forming a blue edge.
		
		Therefore there are at least $32r\sqrt{rc}\cdot n^2$ non-red edges in $C'$. Then, either $C'\subseteq C$, in which case the proof is completed, or $C\subseteq C'$, in which case there are at least 
			\[32r\sqrt{rc}\cdot n^2 - \bigg(\binom{\lceil n/2\rceil}{2} - \binom{|C|}{2}\bigg)\ge (32r\sqrt{rc} - \zeta)n^2\ge 16r\sqrt{rc}\cdot n^2\]
			non-red edges in $C$, as desired.
	\end{proof}
	
	We are ready to conclude the proof of Lemma~\ref{lem:not-red-increasing}. 
	For every vertex $u\in C$, denote by $V_u$ the set of vertices $v\in B$ such that $uv$ is a blue edge in $G$.
	On the one hand, given that there are at~most 
	$\alpha (1/r+\eps) \tbinom{n}{2}\le 2\alpha n^2/r$ blue edges with at least one endpoint in $C$, the number of vertices $u\in C$ with $|V_u|\ge |B|/3$
	is at most $2\alpha n^2/(r\cdot |B|/3) \le 7\alpha n/\sqrt{r}$ by \eqref{eq:bounds_B}. Since there are at most $7\alpha n^2/\sqrt{r} \le 8r\sqrt{rc}\cdot n^2$ edges with at least one such endpoint, at least $8r\sqrt{rc}\cdot n^2$ of the non-red edges guaranteed by Claim~\ref{claim:many-not-red} are between vertices sending fewer than $|B|/3$ blue edges towards $B$.
	On the other hand, for every pair of vertices $u,v\in C$ with $|V_u|,|V_v| < |B|/3$, $u$ and $v$ have at least $|B|/3 - 2\eps n\ge n/4r$ common neighbours $w\in B$ such that none of the edges $uw,vw$ is blue but both are in $G$. Moreover, since $3\zeta n\le n/8r$ by our choice of $c \le 1/2^{39}r^5$, we find that at least half of the vertex pairs between these common neighbours of $u$ and $v$ induce blue edges in $G$. Thus, there are at least $(8r\sqrt{rc}\cdot n^2)\cdot \tbinom{\lfloor n/4r\rfloor}{2}/2 \ge c n^4$ non-red-increasing marked 4-cycles in $G$, as desired.
\end{proof}

\subsection{Nearly balanced \texorpdfstring{$r$-edge-colourings of $K_n$}{r-edge-colourings of K\_n}}
We turn to the analysis of $\eps$-balanced $r$-colourings of $K_n$. 
Our next step is to show that, in any nearly balanced $r$-colouring of $K_n$, a uniformly random embedding of any spanning tree with suitably bounded maximum degree contains approximately the same number of edges in each colour.
In fact, the next more general lemma follows from a simple application of the second moment method.

\begin{lemma}\label{lem:balanced_colours}
	Fix an integer $r\ge 1$, $\delta\in (0,1/2)$, $\eps \in (0,\delta^2/400r)$ and a sufficiently large integer $n$. Consider an $\eps$-balanced $r$-colouring of the edges of $K_n$.
	Let $H$ be a graph on $n$ vertices, $m\ge 1/\eps$ edges and maximum degree $\Delta \le \eps m$.
	Then, with probability at least $1 - 80 \eps r/\delta^2$, a uniformly random embedding of $H$ into $K_n$ contains $(1 \pm \delta)m/r$ edges in each colour.
\end{lemma}
\begin{proof}
	Identify $H$ with its uniformly random embedding and, for every $i\in [r]$, let $X_i$ be the random variable counting the number of edges in $H$ with colour $i$; in particular, $\mathbb E[X_i] = (1/r\pm \eps)m$ for every $i\in [r]$.
    We turn to a computation of $\mathrm{Var}(X_i)$. Note that, for every edge $e$ in $H$, there are at most $2\Delta$ edges in $H$ which share a vertex with $e$.
	For any two vertex-disjoint edges $e,f$ in $H$, denote by $p_1 = p_1(i)$ the probability that $e$ receives colour $i$ and by $p_2 = p_2(i)$ the probability that both $e$ and $f$ receive colour $i$.
	Then, $p_1 = 1/r\pm \eps$ and $p_2 \le p_1\cdot p_1\tbinom{n}{2}/\tbinom{n-2}{2}\le p_1^2 + 5/n$. Thus, for every $i\in [r]$, denoting by $I_e$ the indicator random variable of the event that $e$ receives colour $i$, similarly to~\eqref{eq:var}, we find
	\begin{align*}
	\Var\left(X_i \right) = \sum_{e, f} \left(\mathbb E\left[I_e I_f \right] - \mathbb E\left[I_e \right] \mathbb E\left[I_f \right] \right) = \bigg(\sum_{e} 2\Delta \mathbb E[I_e]\bigg) + \frac{5m^2}{n}\le 2\eps m\cdot m\cdot \bigg(\frac{1}{r} + \eps\bigg) + \frac{5m^2}{n}\le \frac{20\eps m^2}{r}.
	\end{align*}
	
	Using that $\tfrac{1+\delta}{1+\eps r}\ge \tfrac{1-\delta}{1-\eps r}$, Chebyshev's inequality shows that, for every $\eps\in (0, 1/2r)$ and $i\in [r]$, 
\begin{align*}
\mathbb P\left(X_i \neq (1\pm \delta)m/r \right) \le \mathbb P\left(|X_i - \mathbb E[X_i]| \ge (\delta-\eps)m/r \right)\le \frac{\Var(X_i)r^2}{(\delta-\eps)^2 m^2} \le \frac{80\eps r}{\delta^2}.
\end{align*}
Taking a union bound over all $i \in [r]$ finishes the proof.
\end{proof}
	
\section{Proof of Theorems \ref{col:discErdos} and \ref{thm:col_disc_1}}
\label{sec:thm1}

In this section, we prove \Cref{thm:col_disc_1} and then deduce \Cref{col:discErdos} from it.
The proof combines a brief combinatorial case analysis, an application of the probabilistic method and a counting argument and goes roughly as follows:
first, we consider a random embedding of the tree $T$, and then, for some fixed $i\in [r]$, we switch some pairs of edges in order to increase the number of edges in colour $i$.

\begin{proof}[Proof of Theorem~\ref{thm:col_disc_1}]
Fix $i\in [r]$ and an $r$-edge-colouring of $K_n$.
 We start by analysing several particular cases. If $m_i\in \{0, \tbinom{n}{2}\}$, the statement is trivial and, if $\Delta = n-1$, the statement holds by a direct application of the probabilistic method. Suppose that $m_i^*\in [1,16n/3)$ and $\Delta \le n-2$. Then, by choosing $C < 9/256$ so that 
\[C(n-1-\Delta)(m_i^*)^2\le n^3,\quad \text{and consequently}\quad \left\lfloor \frac{2m_i}{n}+1\right\rfloor\ge \frac{2m_i+1}{n}\ge \frac{2m_i}{n} + C(n-1-\Delta)\bigg(\frac{m_i^*}{n^2}\bigg)^2,\]
it is sufficient to show that some copy of $T$ contains a number of edges in colour $i$ different from the expected one: indeed, this guarantees that a copy overshooting this expected number exists.
On the one hand, since $T$ is not a star, Lemma~\ref{lem:many switches} ensures that $T$ contains a switchable pair of edges.
On the other hand, the existence of a non-monochromatic 4-cycle with two edges in colour $i$ opposite to each other guarantees that the above pair of edges can be switched in a suitable copy of $T$ (see Figure~\ref{fig:marked-cycles}(i) and (ii)).
If the latter does not hold, then either all edges in colour $i$ form a clique or contain a common vertex.
In the first case, let the clique have size $k\in [n-1]$. 
Then, we must have $(k-1)k = 2m_i$, which ensures a copy of $T$ with at least $k-1 > 2m_i/n$ edges in colour $i$.
In the second case, $2m_i/n < 2$ but a copy of $T$ with at least two edges in colour $i$ exists, as desired.
 
From now on, we assume that $m_i^*\ge 16n/3$; in particular, $32n/3\le \tbinom{n}{2}$ implying $n\ge 22$.
Let $\mathcal S$ be a free family of $\zeta = \zeta(n,\Delta) = \lceil (n-1-\Delta)/40\rceil$ switchable pairs in $T$ provided by Lemma~\ref{lem:many switches}. 
Recall that a switching is \emph{desirable} with respect to a family $\mathcal C$ of 4-cycles in $K_n$ if it replaces edges $uw$, $vz$ with $uz$, $vw$, where $u,v,w,z$ are embedded into a cycle $\gamma\in \mathcal C$ with $u$ and $v$ at distance 2 in $\gamma$.
    We will take $\mathcal C$ to be the set of cycles for which applying the switching replacing $uw$, $vz$ with $uz$, $vw$ increases the number of edges in colour $i$.
    This is defined more precisely below. 
	Consider an embedding $T'$ of $T$ into $K_n$ chosen uniformly at random.
	Let $X$ be the random variable counting the number edges of colour $i$ in $T'$, and $Y$ be the random variable counting the number of pairs in $\cS$ whose switching is desirable with respect to $\mathcal C$.
	By performing all desirable switchings in $\mathcal S$, we get a copy of $T$ in $K_n$ with at least $X+Y$ edges in colour $i$.
	Hence, there exists a copy $T''$ of $T$ in $K_n$ in which the number of edges in colour $i$ is at least
	\[\mathbb E\left[X+Y \right] = \frac{2m_i}{n} + \mathbb E\left[Y \right]. \]
	In the rest of the proof, we bound $\mathbb E[Y]$ from below.
 
	Fix a pair of edges $(e,f)$ in $\mathcal S$.
	Fix an embedding of $T$ in $K_n$, identify $e$ and $f$ with their images in $K_n$ and let $C\subseteq K_n$ be the 4-cycle along which the allowed switching of $e$ and $f$ can be done.
	This switching of $e,f$ increases the number of edges in colour $i$ if one of the following occurs:
	\begin{itemize}
		\item $C$ contains one edge of colour $i$ and both $e$ and $f$ are of colour different to $i$.
        In this case, the switching of $e,f$ along~$C$ increases the number of edges of colour $i$ by 1.
		\item $C$ contains three edges of colour $i$ and one of $e$ and $f$ is of a different colour.
        Again, the switching of $e,f$ along $C$ increases the number of edges of colour $i$ by 1.
		\item $C$ contains two edges of colour $i$ opposite to each other and both $e$ and $f$ are of a different colour.
        This time, the switching of $e,f$ along $C$ increases the number of edges of colour $i$ by 2.
	\end{itemize}
	The three possibilities coincide with the ones depicted in Figure~\ref{fig:marked-cycles} (where colour $i$ is blue and other colours are are replaced by black dashed lines).
	Now, define $\mathcal C$ to be the family of cycles in $K_n$ coloured in one of the three ways described above.
    Then, since there are $3\tbinom{n}{4}$ 4-cycles in $K_n$, 
	\begin{align}
	\label{eq:expecY}
		\mathbb E[Y] = \frac{|\mathcal S|\cdot |\mathcal C|}{2\cdot 3\tbinom{n}{4}}\ge \frac{(n-1-\Delta)\cdot |\mathcal C|}{10 n^4},
	\end{align}
	where the 2 in the denominator of the expression in the middle comes from the fact that two possible switchings can be done along each 4-cycle.
	
	We now bound $\left|\mathcal C \right|$ from below by counting certain $3$-paths.
	We call a $3$-path \emph{alternating} if its middle edge has a different colour from the other two.
	To count all alternating $3$-paths, we write $m^*_i = \min\left(m_i, \binom{n}{2} - m_i \right)$ and continue with analysing three cases.
    We consider now two colours: colour $i$ and all other colours combined to one colour.
    If $m^*_i=m_i$, colour blue all edges in colour $i$, and red all remaining edges, and if $m^*_i = \binom{n}{2} - m_i$ then colour the edges the other way around.
    Note that we have $m^*_i \le \binom{n}{2}/2$ blue edges in each case.

	For the remainder of this proof, we respectively refer to the number of blue and red edges incident to a vertex as its \emph{blue degree} and \emph{red degree}.
 	Denote by $B$ the subset of vertices of blue degree at least $9n/10$ and let $A = V(K_n) \setminus B$.
	
	\paragraph{Case 1.}
	There are more than $m^*_i/2$ blue edges inside $B$.
	Then we have $|B| \le 5n/9$, as otherwise $B$ is incident to at least $n^2/4$ blue edges, contradicting $m^*_i \le \binom{n}{2}/2$.
	Note that the number of blue edges between $B$ and $A$ is at least $|B|\left(9n/10 - |B| \right)$, and by assumption it is also smaller than $m^*_i/2$.
	Hence,
	\[|B|\left(\frac{9}{10}n - |B| \right) < \frac{1}{2}m^*_i \le \frac{1}{8}n^2. \]
	Solving this quadratic inequality we get that either $|B| > (9+\sqrt{31})n/20 > 5n/9$, which is a contradiction, or $|B| < (9 - \sqrt{31})n/20$.
	In the latter case, we get that each vertex in $B$ has more blue neighbours outside of $B$ than inside of it, which is also a contradiction.
	
	\paragraph{Case 2.}
	There are at least $m^*_i/8$ blue edges inside $A$.
	Since vertices in $A$ have red degree at least $n/10$, every such edge is a middle edge in at least $(n/10)(n/10 - 1)$ alternating $3$-paths.
    Moreover, every marked cycle in $\mathcal C$ contains at most two alternating $3$-paths with blue middle edge, implying that
	\begin{align*}
		\left|\mathcal C \right| \ge \frac{1}{2} \cdot \frac{m^*_i}{8} \cdot \frac{n}{10}\left(\frac{n}{10} - 1 \right) \ge \frac{m^*_i n^2}{3200} \ge \frac{(m^*_i)^2}{800},
	\end{align*}
    as $n \geq 20$ and $m^*_i \le \binom{n}{2} / 2$. Plugging this inequality into \eqref{eq:expecY} finishes the proof in this case.
	
	\paragraph{Case 3.} There are at least $3m^*_i/8$ blue edges between $A$ and $B$.
	Then, there must be at least $3m^*_i/8|A| \ge 3m^*_i/8n \ge 2$ vertices in $B$.
	Let $u,v$ be two distinct vertices in $B$ and let $J$ be the intersection of their neighbourhoods. By definition of $B$, we have that
 $|J| \ge \left(2\cdot9/10-1 \right)n = 4n/5$.
	Moreover, the number of red edges inside $J$ is at least $\binom{|J|}{2} - m^*_i$, and every such red edge is a middle edge in a alternating $3$-path with $u,v$ as endpoints.
	In total, we get that the number of such alternating $3$-paths is at least
	\[\binom{3m^*_i/8n}{2} \cdot \left(\binom{4n/5}{2} - m^*_i \right). \]
	Again, every marked cycle in $\mathcal C$ contains at most two such alternating $3$-paths, implying that
    \begin{align}
    \label{eq:ineqCcase3}
		\left|\mathcal C \right| &\ge \frac{1}{2} \binom{3m^*_i/8n}{2} \cdot \left(\binom{4n/5}{2} - m^*_i \right).
    \end{align}
    Combining the bounds $\binom{3m^*_i/8n}{2} \geq \frac{(3m^*_i/8n)^2}{4}$ (using that $m^*_i \geq 16n/3$), $\binom{4n/5}{2} \geq \frac{4n^2}{15}$ (using that $n \geq 8$), $m^*_i \leq n^2/4$ and \eqref{eq:ineqCcase3} gives
    \begin{align*}
        \left|\mathcal C \right| \ge \frac{1}{2} \cdot \frac{1}{4} \cdot \frac{9(m^*_i)^2}{64n^2} \left (\frac{4n^2}{15} -\frac{n^2}{4} \right ) \ge \frac{(m^*_i)^2}{10^4}.
    \end{align*}
    In total, by plugging in those bounds in~\eqref{eq:expecY} we get
	\[\mathbb E\left[Y \right] \ge \frac{1}{10^5}\left(n - 1 - \Delta \right) \left(\frac{m^*_i}{n^2} \right)^2 , \]
	finishing the proof.
\end{proof}

We finish this section by deducing \Cref{col:discErdos} from \Cref{thm:col_disc_1}.

\begin{proof}[Proof of \Cref{col:discErdos}.]
Fix $r\ge 2$ and an edge-colouring of $K_n$ in $r$ colours. 
Let $i\in [r]$ be one of the most represented colours; in particular, $m_i \geq \frac{1}{r}\binom{n}{2}$.
If $m_i\ge \frac{1.5}{r}\binom{n}{2}$, by choosing $c$ suitably small, the statement follows from the probabilistic method.
Otherwise, $m_i\in [\frac{1}{r}\binom{n}{2}, \frac{1.5}{r}\binom{n}{2}]$.
Then, by \Cref{thm:col_disc_1}, there exists a copy of $T$ in which the number of edges in colour $i$ is at least
\[\frac{2m_i}{n} + C\left(n - 1 - \Delta \right)\left(\frac{m_i^*}{n^2} \right)^2 \ge \frac{n-1}{r} + C \min\bigg(\frac{1}{r^2}, \bigg(1-\frac{1.5}{r}\bigg)^2\bigg)\left(n - 1 - \Delta \right).\]
Choosing $c = C/4$ finishes the proof.
\end{proof}

\section{Proof of Propositions~\ref{prop:example_many_colours} and~\ref{prop:k>2}}\label{sec:3}

Fix $\eps > 0$. Before we proceed with the proof of Proposition~\ref{prop:example_many_colours}, we recall that, given a prime number $p$, 
there is a finite projective plane $\mathrm{PG}(2,p)$ with $p^2+p+1$ points and $p^2+p+1$ lines where every line contains $p+1$ points, every point lies on $p+1$ lines and every two lines intersect at a single point.
We use this classical geometric object in combination with the following result for the density of primes.

\begin{lemma}[Theorem 1 in~\cite{BHP01}]\label{lem:gaps_primes}
	Let $(p_n)_{n\ge 1}$ be the sequence of prime numbers listed in increasing order. Then, for all sufficiently large $n$, $p_{n+1} - p_n \le p_n^{0.525}$.
\end{lemma}

We are now ready to prove Proposition~\ref{prop:example_many_colours}.

\begin{proof}[Proof of Proposition~\ref{prop:example_many_colours}]
	Let $\alpha = 0.525$ and $\eps$ be sufficiently small so that, by \Cref{lem:gaps_primes}, we have $p_{m+1}-p_m\le p_m^{\alpha}$ for every prime $p_m$ satisfying $p_m\ge (80/\eps)^{1/(1-\alpha)}$. 
    Let $p_{m_0}$ be the smallest prime greater than or equal to $(80/\eps)^{1/(1-\alpha)}$.
	Set $r_0 = p_{m_0}^2+p_{m_0}+2$ and fix any $r\ge r_0$. Then, there is a prime number $p\ge p_{m_0}$ such that 
	\begin{equation*}
	p^2+p+2\le r\le (p+p^{\alpha})^2 + (p+p^{\alpha}) + 2\le (1+3p^{\alpha-1})(p^2+p+2).
	\end{equation*}
	Setting $k = p^2+p+2$, the above rewrites as
    \begin{equation}\label{eq:bounds_r_minus_k}
        0 \le r-k \le 3p^{\alpha-1} k.
    \end{equation}
 
    We first define a suitable $k$-colouring of $K_n$ and then modify it to form an $r$-colouring.
	To begin with, we identify the set of colours $\{2,\ldots,k\}$ with the set of points in $\mathrm{PG}(2,p)$ via an arbitrary bijection.
	Also, we assume for convenience that $n$ is divisible by $(r-1)(k-1)$; the general case follows by putting lower and upper integer parts wherever is necessary.
	
	Now, we describe the $k$-colouring $\mathbf{c}_{k}$ of $K_n$ we are going to work with. 
	We begin by partitioning the vertices of $K_n$ into $k-1$ sets $V_2, \ldots, V_{k}$ of equal size (called \emph{parts}). 
	First, we give colour 1 to all edges whose both endpoints are in the same part.
	Second, we identify the parts with the lines in $\mathrm{PG}(2,p)$ in a bijective manner and, for every pair of distinct integers $i,j\in [2,k]$, we colour all edges between $V_i$ and $V_j$ in the colour associated to the unique point in $\mathrm{PG}(2,p)$ where the lines corresponding to $V_i$ and $V_j$ intersect.
	In this way, every vertex is adjacent to edges in exactly $p+2$ colours.
	
	Next, we modify the $k$-colouring $\mathbf{c}_{k}$ to create a suitable $r$-colouring $\mathbf{c}_r$ of $K_n$. 
	For every pair of distinct integers $i,j\in [2,k]$, we partition the edges between $V_i$ and $V_j$ into $n/(k-1)$ perfect matchings.
	Then, we further partition this set of matchings into $r-1$ groups of equal size and, for every $i'\in [k+1,r]$, recolour all (edges in) matchings in the $i'$-th group in colour $i'$. 
	By construction, every colour in $[2,r]$ is given to
    \begin{align} \label{eq:NumberColours2tor}
        \frac{1}{r-1}\left(\binom{n}{2} - (k-1)\binom{n/(k-1)}{2}\right) = \frac{1}{r-1} \left(1 - \frac{n/(k-1)-1}{n-1}\right) \binom{n}{2}\le \frac{1}{r} \binom{n}{2}
    \end{align}
	edges in $\mathbf{c}_r$, where the difference between the last two expressions above is of order $\Theta(n^2)$. 
	In particular, in $\mathbf{c}_r$, all colours in $[2,r]$ are given to the same number of edges while colour 1 is given to $\Theta(n^2)$ more edges than the rest.
	
	Fix $L \defined \lceil 1/\eps\rceil\le 1+1/\eps$, an $L$-caterpillar $T$ and 
	let $s=\lceil (L+1)/n \rceil$.
	Then, every vertex in the central path of $T$ has degree $s$ or $s+1$.
	Fix an arbitrary copy of $T$ in $K_n$. 
	On the one hand, it can contain at most $(L+1)\cdot (n/(k-1)-1)$ edges in colour 1 in $\mathbf{c}_r$: indeed, every such edge is incident to a vertex $v$ in the central path of $T$ and its other endvertex must therefore be in the same part as $v$.
	On the other hand, every vertex $u$ in $K_n$ is incident to $n/(k-1)-1$ edges in colour 1, and $(k-2)n/(k-1)(r-1)$ edges in each colour among $k+1,\ldots,r$ (by rearranging the formula given by \eqref{eq:NumberColours2tor}) while all remaining edges incident to $u$ have one of $p+1$ colours among $2,\ldots,k$. 
	Moreover, it follows from $p \geq (80/\eps)^{1/(1-\alpha)}$ and \eqref{eq:bounds_r_minus_k} that $\max(2(p+1)(L+1)(L+2), 2(L+2)(r-k+1))\le k-1$.
    Combining this with $s\cdot (L+1)\ge n$, we get
	\begin{align*}
	\frac{1}{p+1}\left(s - \left(\frac{n}{k-1}-1\right) - \frac{(r-k)(k-2)n}{(k-1)(r-1)}\right)
	&> \frac{1}{p+1} \left(s - \frac{n}{k-1} - \frac{(r-k)n}{(k-1)}\right)\\
	&\ge \frac{s}{2(p+1)}\ge (L+1) \cdot \left(\frac{n}{k-1}-1\right).    
	\end{align*}
	We deduce that, in $\mathbf{c}_r$, some colour among $2,\ldots,k$ is given to more edges in the current copy of $T$ than colour $1$. 
	Finally, to account for the fact that $\mathbf{c}_r$ is not a balanced colouring (colour 1 is given to more edges than the rest), it suffices to reassign the colour of some edges with $\mathbf{c}_r$-colour 1. 
	Recall that $\eps$ has been chosen sufficiently small so that $p_{m+1}-p_m\le p_m^{\alpha}\le p_m$ is satisfied for every $p_m\ge (80/\eps)^{1/(1-\alpha)}$ and, therefore, 
	\[r_0 = p_{m_0}^2+p_{m_0}+2\le (2(80/\eps)^{1/(1-\alpha)})^2 + 2(80/\eps)^{1/(1-\alpha)} + 2 \le 1/\eps^5,\] 
	as desired.
\end{proof}

Now, we transition to the proof of Proposition~\ref{prop:k>2}.

\begin{proof}[Proof of Proposition~\ref{prop:k>2}]
	We assume for convenience that $n$ is divisible by $2(r-1)$ and show the statement for $C=1$; the general case requires only minor technical modifications resulting in a larger constant $C$.
	Let $V_2,\ldots,V_r$ be a partition of the vertex set of $K_n$ into $r-1$ parts of size $t = n/(r-1)$. 
	We consider two cases.
	
	We define a balanced $r$-colouring $\mathbf{c}_r$ of $K_n$ as follows. 
	To begin with, for every pair of $i,j\in [2,r]$ with $i<j$, decompose the edges between $V_i$ and $V_j$ into $t$ matchings $M_1, \ldots, M_t$ of size $t$. 
	Then, give colour $i$ to the edges in $M_1\cup \ldots \cup M_{t/2}$ and colour $j$ to the edges in $M_{t/2+1}\cup \ldots \cup M_t$.
	In this way, every colour in $[2,r]$ is given to $(r-2)t^2/2$ edges.
	Finally, for every $i\in [2,r]$, define $d = \lfloor(n-(r-1)^2)/r(r-1)\rfloor$,
	pick an arbitrary $d$-regular subgraph of $K_n[V_i]$ and give colour $i$ to each of its edges while all remaining edges in $K_n$ are given colour 1.
	It is easy to check that our choice of $d$ satisfies
	\[\frac{(r-2)t^2}{2} + \frac{dt}{2}\le \frac{1}{r}\binom{n}{2},\]
	so colours $2,\ldots,r$ are given to the same number of edges (equal to the left-most term in the above display) while colour 1 is given to the same or a slightly larger number of edges.
	
	Now, fix a tree $T$ of maximum degree $\Delta\ge 23 n/24 + 1$ and a copy of $T$ in $K_n$. 
	On the one hand, the number of edges in $T$ with colour 1 in $\mathbf{c}_r$ is at most 
	\begin{equation}\label{eq:1}
	t - 1 - d + (n-1-\Delta).
	\end{equation}
	On the other hand, every vertex $v$ of degree $\Delta$ in our copy of $T$ is incident to at least
	\begin{equation}\label{eq:2}
	\Delta - \left(t-1-d\right) - (r-2)t/2
	\end{equation}
	edges in the same colour (given by the index of the part containing $v$).
	A direct computation shows that $\eqref{eq:1} < \eqref{eq:2}$ when $\Delta\ge (3r^2-4)n/4r(r-1)+1$. 
	Moreover, the function $x\in [3,\infty)\mapsto (3x^2-4)/4x(x-1)$ is decreasing and its value at $x=3$ is $23/24$, so the condition on $\Delta$ from the statement of the lemma is sufficient to guarantee that, for any copy of $T$, some colour among $2,\ldots, r$ is given to more edges than colour 1.
    Finally, to account for the fact that $\mathbf{c}_r$ is not a balanced colouring (colour 1 is given to more edges than the rest), it suffices to reassign the colour of some edges with $\mathbf{c}_r$-colour 1. This concludes the proof.
\end{proof}

\section{Proofs of Theorems \ref{thm:col_disc_2} and \ref{thm:col_disc_3}}
\label{sec:5}

As in the proof of Theorem~\ref{thm:col_disc_1}, the existence of many pairs of edges in $T$ that can be switched independently plays a key role in the proof of Theorems~\ref{thm:col_disc_2} and~\ref{thm:col_disc_3}.
However, this time, we need to control the number of edges in \emph{each} colour, not only a single one. 
We start with the proof of \Cref{thm:col_disc_2}.

\begin{proof}[Proof of Theorem~\ref{thm:col_disc_2}]
Fix $r\ge 3$ and $\eta\in (0,2^{-15}/r)$.
For every $i\in [2,r]$, consider the family $\mathcal C_i$ of non-$i$-increasing switchings introduced in Section~\ref{sec:4.2} (see Figure~\ref{fig:marked-cycles}),
and define $c = \eta^2/2^{10}r^3$, $\delta = c/200\ge \eta^2/10^6 r^3$ and $\eps =\delta^2/1000r$, so that the assumptions of Lemmas~\ref{lem:marked_4_cycles},~\ref{lem:not-red-increasing} and~\ref{lem:balanced_colours} are simultaneously satisfied, as well as 
$400\eps r/c+80\eps r/\delta^2 < 1$.
Then, there is a copy of $T$ in $K_n$ where each colour is given to $(1\pm \delta)(n-1)/r$ edges and moreover, there is a free family containing at least $cn/40$ non-$i$-increasing switchings for every $i\in [2,r]$.
After doing each of these switchings, the difference between number of edges in colour 1 and the number of edges in colour $i$ becomes at least $cn/40 - 2\delta n/r\ge \delta n$ for every $i\in [2,r]$, which finishes the proof.
\end{proof}

\subsection{Dealing with the missing edges: deducing Theorem~\ref{thm:col_disc_3}}

To prove Theorem~\ref{thm:col_disc_3}, we need to extend our techniques to sparser host graphs. 
Given a graph $G$, we say that a vertex pair $uv$ is a \emph{ghost edge (for $G$)} if $uv$ is not an edge of $G$.
We will perform an additional round of vertex exchanges to ensure that there are no ghost edges remaining in our final copy of $T$.

When we say that two copies $T',T''$ of $T$ \emph{differ in at most $\ell$ edges}, we mean that the implicit isomorphisms $f\from T\to T'\subseteq K_n$ and $g\from T\to T''\subseteq K_n$ have different images for at most $\ell$ edges of $T$. 
We also say that $T''$ is obtained from $T'$ by \emph{exchanging the vertices $u,v\in V(K_n)$} if 
$f(g^{-1}(u))=v$, $f(g^{-1}(v)) = u$ and, moreover, $f$ and $g$ coincide over $V(T)\setminus f^{-1}(\{u,v\}) = V(T)\setminus g^{-1}(\{u,v\})$.

\begin{lemma}
	\label{lem:remove-ghost-edges}
	Fix $\eps\in (0,1/2(\Delta+1)^3)$, $\delta > 0$, a graph $G\subseteq K_n$ with $n$ vertices and minimum degree $\delta(G)\ge (1-\eps)n$, and a copy $T'\subseteq K_n$ of a graph $T$ containing at most $\delta n$ ghost edges.
	Then, there is a copy $T''\subseteq G$ of $T$ such that $T'$ and $T''$ differ in at most $2\Delta\delta n$ edges.
\end{lemma}
\begin{proof}
	If $n < 2(\Delta+1)^3$, there is nothing to prove (since $G=K_n$), so suppose that $n\ge 2(\Delta+1)^3$.
	We construct a sequence of copies $T' = T_0, T_1,\ldots, T_k$ of $T$ in $K_n$ such that, for every $i\in [k]$, $T_i$ contains fewer ghost edges than $T_{i-1}$ and $T_i$ differs from $T_{i-1}$ in at most $2\Delta$ edges. 
	Moreover, we will ensure that $T_k$ contains no ghost edges, which will allow us to choose $T'' = T_k$: note that, by construction, $k\le \delta n$, so $T'$ and $T''$ will indeed differ in at most $2\Delta\delta n$ edges.
	
	We only describe the construction of $T_1$ from $T_0$, as the subsequent constructions are done similarly.
    Let $N_H[x]$ denote the \emph{closed $H$-neighbourhood} of $x$, which consists of $x$ and all its neighbours in $H$.
	Suppose that $x$ is incident to a ghost edge in $T_0$. 
	We will show that there is a vertex $z$ at distance at least 3 from $x$ in $T_0$ such that none of the edges between $x$ and $N_{T_0}[z]$ is a ghost edge, and none of the edges between $z$ and $N_{T_0}[x]$ is a ghost edge. 
    We then obtain the tree $T_1$ by exchanging the vertices $x$ and $z$ in $T_0$; clearly, in this process, at least one ghost edge is eliminated.
	
	To find a vertex $z$ as above, we first construct a set $S$ of vertices containing $x$ and of size at least $\lceil n/(1+\Delta+\Delta^2)\rceil$, such that any two vertices in $S$ are at distance at least three in $T_0$. 
	We do this by performing the following greedy deletion procedure. 
	At the first stage, we add $x$ to $S$ and delete all vertices within distance 2 from $x$. 
	Then, as long as there is at least one not yet deleted vertex, we choose an arbitrary such vertex $y$, add it to $S$ and delete all remaining vertices within distance 2 from $y$ in $T_0$.

    Observe that for every distinct $z,z' \in S$, we have that $N_{T_0}[z]$ and $N_{T_0}[z']$ are disjoint.
    Additionally, there are at most $(\Delta+1)\eps n < n/(1+\Delta+\Delta^2) - 1$ ghost edges going out of $N_{T_0}[x]$.
    As a result, there is a vertex $z$ in $S$ such that the bipartite graph with parts $N_{T_0}[x]$ and $N_{T_0}[z]$ does not contain any ghost edge, which is therefore a complete bipartite graph.
	This choice of $z$ satisfies the desired properties. 
\end{proof}

\begin{remark}\label{rem:non-tree}
		Note that the statement of Lemma~\ref{lem:remove-ghost-edges} is valid for any graph $T$ of maximum degree $\Delta$ (in particular, ones containing cycles). Moreover, we note that, at the last step of the proof of Lemma~\ref{lem:remove-ghost-edges}, there are at least 
		\[\frac{n}{1+\Delta+\Delta^2} - (\Delta+1)\eps n\ge \frac{n}{1+\Delta+\Delta^2} - \frac{n}{2(\Delta+1)^2}\ge \frac{n}{2(\Delta+1)^2}\]
		choices for $z$. This observation will be useful in Section~\ref{sec:5'}.
\end{remark}

We now apply \Cref{lem:remove-ghost-edges} to prove \Cref{thm:col_disc_3}.

\begin{proof}[Proof of \Cref{thm:col_disc_3}]
	Fix $r\ge 3$, $\Delta\ge 2$ and $\eta\in (0,2^{-15}/r)$, and define $c = \eta^2/2^{10}r^3$, $\eps = c^2/10^6 r$ and $\delta = 30\sqrt{r\eps}$. In particular, we have that
 \[\frac{80\eps r}{\delta^2} + \frac{400\eps r}{c} < 1\quad \text{and}\quad 4\Delta \eps + 20r\eps \le \frac{c}{160}.\]
	Moreover, fix an $\eps$-balanced $\eta$-standard $r$-colouring of $G$ and extend it arbitrarily to an $\eps$-balanced $\eta$-standard $r$-colouring of $K_n$. 
	For every $i\in [2,r]$, denote by $\mathcal{C}_i'$ the family of marked 4-cycles in $G$ which are not $i$-increasing.
	Recall that, for every $i\in [2,r]$, Lemma~\ref{lem:not-red-increasing} implies that $|\mathcal{C}_i'|\ge cn^4$.
	
	Consider a uniformly random copy $T'$ of $T$ in $K_n$ and let $X$ be the number of ghost edges in $T'$. 
	Since $\mathbb E[X]\le \eps n$, a computation similar to~\eqref{eq:var} shows that $\mathrm{Var}(X)\le 2\Delta n + O(n)$, which implies that $X\le 2\eps n$ with high probability.
	Together with Lemma~\ref{lem:marked_4_cycles} applied for $\mathcal{C}_2', \ldots, \mathcal{C}_r'$ and Lemma~\ref{lem:balanced_colours}, we deduce that, with positive probability, 
	$T'$ contains at most $2\eps n$ ghost edges, between $(1-\delta)(n-1)/r$ and $(1+\delta)(n-1)/r$ edges in each colour, and a free family with at least $cn/40$ pairs in non-$i$-increasing switchings for all $i\in [2,r]$. 
	We condition on these properties of $T'$.
	
	Using Lemma~\ref{lem:remove-ghost-edges}, one can produce a copy $T''$ of $T$ in $K_n$ that contains no ghost edges and differs from $T'$ in at most $2\Delta \eps n$ edges (and thus, the symmetric difference of $E(T')$ and $E(T'')$ contains no more than $4\Delta \eps n$ edges). 
	As a result, since $4\Delta \eps < c/80$, one can find a free family $\mathcal S$ containing at least $cn/80$ pairs in non-$i$-increasing switchings in $T''\subseteq G$ for every $i\in [2,r]$. 
	Moreover, for every $i\in [2,r]$, $T'$ contains at most $2\delta n$ more edges in colour $i$ than in colour 1, which means that after performing all switchings in $\mathcal S$, the resulting copy of $T$ contains at least $(c/40 - 4\Delta \eps - 2\delta)n\ge cn/160 \ge\eps n$ more edges in colour 1 than in colour~$i$, as desired.
\end{proof}

\section{\texorpdfstring{$k$}{k}-sum subgraphs and inert triplets: proofs of Proposition~\ref{prop:par-bip} and Theorem~\ref{thm:AlmostZeroSumDenseGraph}}\label{sec:5'}

We start with a short proof of Proposition~\ref{prop:par-bip}.

\begin{proof}[Proof of Proposition~\ref{prop:par-bip}]
	Fix an inert triplet $(H,G,\bc)$ and extend the bipartite colouring $\bc$ to a complete bipartite colouring of the edges of $K_n$ with parts $V_1$ and $V_2$.
	Note that, for every copy $H'$ of $H$ in $K_n$, the number of vertices of odd degree in $H'[V_1]$ is even.
	As a result, the parity of the number of edges of $H'$ in the cut $(V_1,V_2)$ is fixed, which finishes the proof since $\bc(H') = |E(H)| - 2|E(H')\cap (V_1\times V_2)|$. 
\end{proof}

Proposition~\ref{prop:par-bip} and the simple fact that $\bc(H')$ has the same parity for all copies $H'$ of $H$ in $G$ show the ``only if'' part of Theorem~\ref{thm:AlmostZeroSumDenseGraph}. 
We turn to the more interesting ``if'' part of the theorem.
If the triplet $(H,G,\bc)$ is not inert, our first task is to find a small subgraph of the graph $G$ where one vertex exchange can change the remainder modulo 4 of $\bc(H')$ for a suitable copy $H'$ of $H$.
Given a $2$-edge-colouring $\bc$ of a graph $G$ in red and blue, call a $4$-cycle in $G$ a \emph{$(1,3)$-cycle} if it has three edges in one colour and one edge in the other. 
We start with the following auxiliary structural lemma.

\begin{lemma}\label{lem:Ramsey}
	Let $\eps\in [0,1/10]$, and consider a sufficiently large $n$ and an $n$-vertex graph $G$ with minimum degree at least $(1-\eps)n$, along with a non-bipartite $2$-edge-colouring $\bc$ of $G$. Then, $G$ contains a $(1,3)$-cycle.
\end{lemma}

\begin{proof}
	Suppose for contradiction that no $(1,3)$-cycle exists.
	To begin with, the classical Tur\'an's theorem implies that $G$ contains a copy of the complete graph $K_6$.
	Moreover, since every 2-edge-colouring of $K_6$ contains a monochromatic triangle, there is a monochromatic triangle $T$ in $G$.
	Then, to avoid a $(1,3)$-cycle, every vertex $v$ in $G\setminus T$ must have all of its edges towards $T$ in the same colour.
	
	Suppose without loss of generality that all edges in $G[T]$ are red. 
	Denote by $R$ (resp.\ $B$) the set of vertices in $G\setminus T$ with three red (resp.\ blue) edges towards $T$, and set $X = V(G)\setminus (R\cup B\cup T)$.
	Then, since no $(1,3)$-cycles exist, all edges in $G[R]\cup G[B]$ must be red.
	Furthermore, since the minimum-degree condition on $G$ implies that $|X|\le 3\eps n$, there exists an inclusion-maximal set $Y$ of size at least $(n-3-3\eps n)/2\ge n/3$ such that all edges in $G[Y]$ are red.
	
	We claim that, for every vertex $x$ in $G\setminus Y$, all edges between $x$ and $Y$ are blue.
	Suppose otherwise for contradiction. 
	Then, by maximality of $Y$, there must be vertices $x$ in $G\setminus Y$ and $y,z\in Y$ such that $xy, xz$ are respectively red and blue.
	Since $|Y|\ge 2\eps n+3$, the vertices $y$ and $z$ must have a common neighbour $w$ in $Y$. 
	However, the cycle $xywz$ would then form a $(1,3)$-cycle, leading to a contradiction.
	Thus, all edges between $G\setminus Y$ and $Y$ are blue.
	
	Finally, we claim that all edges within $G\setminus Y$ are red. 
	Suppose for contradiction that there is a blue edge $xy$ in $G\setminus Y$. Let $z$ be a neighbour of $x$ in $Y$. 
    Since $y$ and $z$ have a common neighbour $w$ in $Y$, the cycle $xywz$ would then form a $(1,3)$-cycle, contradicting the assumption.
	Hence, the colouring of $G$ is bipartite with parts $Y$ and $V(G)\setminus Y$, contradicting the assumptions of the lemma and confirming the existence of a $(1,3)$-cycle.
\end{proof}

Consider a graph $G$ with a $2$-edge-colouring $\bc$, an integer $\Delta\ge 1$ and a copy $K$ of the complete bipartite graph $K_{2,2\Delta}$ in $G$. 
Suppose that $K$ has parts $A = \{x,z\}$ and $B = \{y\}\cup B'$ with $|B| = |A|\Delta = 2\Delta$.
We call $K$ a \emph{type-$1$} copy (of $K_{2,2\Delta}$) if:
\begin{enumerate}[(i)]
	\item\label{pt:i} each of the paths $\{xwz: w\in B'\}$ is coloured identically, and
	\item\label{pt:ii} $\bc(xy)\neq \bc(yz)$ and the path $xyz$ is coloured differently from the other $2\Delta-1$ paths of length 2 with endpoints $x,z$, that is, $\bc(xy)\neq \bc(xw)$ or $\bc(yz)\neq \bc(wz)$ for all $w\in B'$.
\end{enumerate}
Moreover, we call $K$ a \emph{type-$2$} copy (of $K_{2,2\Delta}$) if~\eqref{pt:i} holds and, additionally, for all $w\in B'$, $|\bc(\{xy,yz\})|\neq |\bc(\{xw,wz\})|$.
In particular, in a type-$2$ copy, the path $xyz$ is monochromatic if and only if the paths $xwz$ with $w\in B'$ are not.

Our gadgets (resp.\ 2-gadgets) mentioned in the proof outline will consist of parts of the graph $H$ embedded into a type-1 (resp.\ type-2) copy of $K_{2,2\Delta}$ in $G$. The next lemma shows that a type-2 copy of $K_{2,2\Delta}$ can be found in dense host graphs equipped with non-bipartite colourings.

\begin{lemma}\label{lem:2-gadget}
	Fix $\eps$, $n$, $G$ and $\bc$ as in Lemma~\ref{lem:Ramsey}. 
	Then, $G$ contains a type-$2$ copy of $K_{2,2\Delta}$.
\end{lemma}
\begin{proof}
	Using \Cref{lem:Ramsey}, let $xyzw$ be a $(1,3)$-cycle where the edge $xy$ has a different colour than the other edges of the cycle. Then, since at least $(1-2\eps)n - 5\ge 8\Delta-5$ vertices different from $x,y,z,w$ are adjacent to both $x$ and $z$, there must be a set $B'$ of size $2\Delta-1$ such that the paths $\{xwz:w\in B'\}$ are coloured identically. 
	If these paths are monochromatic, the edges between $\{x,z\}$ and $\{y\}\cup B'$ form a type-2 copy of $K_{2,2\Delta}$. Otherwise, the edges between $\{x,z\}$ and $\{w\}\cup B'$ form a type-2 copy of $K_{2,2\Delta}$.
\end{proof}

We now formally define our gadgets.

\begin{definition}\label{gadgets}
	In the setting of Theorem~\ref{thm:AlmostZeroSumDenseGraph}, given a copy $H'$ of $H$ in $G$, we define a \emph{gadget} (resp.\ \emph{$2$-gadget}) to be a type-1 (resp.\ type-2) copy $K$ of $K_{2,2\Delta}$ in $G$ with parts $A = \{x,z\}$ and $B = \{y\}\cup B_1\cup B_2\cup B'$ for disjoint sets $B_1,B_2,B'$ such that
	\[N_{H'}(x) = \{y\}\cup B_1\quad \text{and}\quad N_{H'}(z) = B_2\]
	have the same size.
\end{definition}

We highlight the following observation.

\begin{remark}
\label{rem:PositiveNegativeGadget}
In a gadget (resp.\ a $2$-gadget), exchanging the positions of the vertices $x,z$ in the embedding $H'$ of $H$ produces an embedding $H''$ of $H$ in $G$ with $\bc(H'') = \bc(H')\pm 4$ (resp.\ $\bc(H'') = \bc(H')\pm 2$) but $\bc(H'')\neq \bc(H')$. 
\end{remark}
\noindent
Additionally, a gadget is called \emph{positive} if $\bc(H'') > \bc(H')$ and \emph{negative} otherwise.

Our next goal is to confirm that dense host graphs contain many type-1 copies of $K_{2,2\Delta}$. The next lemma is a preliminary step towards the proof of this statement.

\begin{lemma}
	\label{lem:linearly-balanced-vertices}
	Fix $\eps, \nu > 0$, and consider an $n$-vertex graph $G$ with minimum degree at least $(1-\eps)n$ and a $\nu$-balanced $2$-edge-colouring $\bc$ of $G$. 
	Set $\xi := (1 - 2\nu - 4\eps)/8$. 
	Then, in $G$, there are at least $\xi n$ vertices, each having at least $\xi n$ blue neighbours and at least $\xi n$ red neighbours in $G$.
\end{lemma}
\begin{proof}
	Assume for contradiction that the statement does not hold.
	Consider a partition $R\union B\union Y$ of $V(G)$ where vertices in $Y$ are incident to at least $\xi n$ edges in each colour, vertices in $R$ are incident to fewer than $\xi n$ blue edges and vertices in $B$ are incident to fewer than $\xi n$ red edges.
	Suppose without loss of generality that $|R|\le |B|$. 
    We consider two cases.
	If $|R|\ge (2\xi + \eps)n$, since every vertex in $R$ (resp.\ $B$) is incident to fewer than $\xi n$ blue (resp.\ red) edges, there are fewer than $(|R|+|B|)\xi n\le 2|B|\xi n$ edges between $R$ and $B$. However, this implies that some vertex in $B$ has more than $\eps n$ non-neighbours in $R$, contradicting the minimum-degree condition for $G$.
	If $|R| < (2\xi + \eps) n$, then $|B|\ge n-|R\cup Y|\ge (1-3\xi-\eps)n$ and each vertex in $B$ is incident to at least $(1-\eps)n$ edges, of which less than $\xi n$ are red, and so at least $(1 - \eps - \xi)n$ are blue. Hence, the number of blue edges in $G$ is at least
	\[\frac{(1-3\xi-\eps)n\cdot (1-\xi-\eps)n}{2} > (1-4\xi -2\eps)\binom{n}{2} \ge \bigg(\frac{1}{2}+\nu\bigg)\binom{n}{2},\]
which contradicts the assumption that $\bc$ is a $\nu$-balanced colouring. This completes the proof.
\end{proof}

\begin{lemma}
	\label{lem:lots-of-gadget-sites}
	Fix $\eps,\nu\in [0,1/800]$, $\Delta\ge 1$ and consider an $n$-vertex graph $G$ with minimum degree at least $(1-\eps)n$ and a $\nu$-balanced $2$-edge-colouring $\bc$ of $G$. 
    Then, $G$ contains at least $n^{2\Delta+2}/(40\Delta)!$ type-$1$ copies of $K_{2,2\Delta}$.
\end{lemma}
\begin{proof}
	For a vertex $u$, let $N_B(u)$ (resp.\ $N_R(u)$) denote the set of neighbours of $u$ that are connected to $u$ via a blue (resp.\ red) edge in $G$.
	By \Cref{lem:linearly-balanced-vertices}, there are at least $n/9$ vertices $x$ such that $\min(|N_B(x)|, |N_R(x)|)\ge n/9$.
	Fix one such vertex $x$ and assume without loss of generality that at least half of the edges between $N_B(x)$ and $N_R(x)$ are blue.
	Consequently, there must be at least $(n/9)(n/9-n/800)/2\ge n^2/170$ blue edges between $N_B(x)$ and $N_R(x)$.
    Therefore there must be at least $n/340$ vertices $y\in N_R(x)$ such that $|N_B(x)\cap N_B(y)|\ge n/340$.
	
	Fix one such vertex $y\in N_R(x)$ and let $Y = N_B(x)\cap N_B(y)$.
	Also, fix a vertex $z\in Y$ and suppose that $|N_B(z)\cap Y|\ge |N_R(z)\cap Y|$ (the opposite case is treated similarly).
	Then, we have $|N_B(z)\cap Y|\ge (n/340-\eps n-1)/2\ge n/1600$ and, for every subset $W\subseteq N_B(z)\cap Y$ of size $2\Delta - 1$, the complete bipartite graph between $\{x,z\}$ and $\{y\}\cup W$ is a type-1 copy of $K_{2,2\Delta}$. 
	As a result, the number of such copies is at least
	\[\frac{n}{9}\cdot \bigg(\frac{n}{340}\bigg)^2\cdot \binom{ n/1600 }{2\Delta-1}\ge \frac{n^{2\Delta+2}}{2^{11\cdot (2\Delta+2)}(2\Delta)!}\ge \frac{n^{2\Delta+2}}{(40\Delta)!},\]
as desired.
\end{proof}

For the remainder of this section, given a host graph $G\subseteq K_n$ with a $\nu$-balanced $2$-edge-colouring $\bc:E(G)\to \{-1,1\}$, we consider an extension of $\bc$ by defining $\bc':E(K_n)\to \{-1,1\}$ to be a $\nu$-balanced $2$-edge colouring of $K_n$.
Our next goal is to show that, given an $n$-vertex graph $G$ with minimum degree at least $(1-\eps)n$ and a graph $H$ on $m$ edges and maximum degree $\Delta$, with positive probability, a uniformly random embedding $H'$ of $H$ in $K_n$ simultaneously satisfies the following properties: 
\begin{enumerate}
	\item[(a)] $H'$ has $(1\pm \nu/5)m/2$ edges in each colour,
	\item[(b)] $H'$ contains at most $3\eps |E(H)|$ edges missing from $G$, and
	\item[(c)] $H'$ contains at least $m/(45\Delta)!$ positive and $m/(45\Delta)!$ negative vertex-disjoint gadgets.
\end{enumerate}
Note that property (a) holds with probability at least $0.8$ by Lemma~\ref{lem:balanced_colours} applied with $\delta = \nu/5\in (0,2^{-20}]$ and $\eps\le \nu^2/10^5$, and property (b) holds with probability at least $2/3$ by Markov's inequality.
Our next lemma verifies that property (c) also holds with suitable probability.

\begin{lemma}
	\label{lem:find-common-gadgets}
	Fix $\eps,\nu\in [0,1/800]$, $m\ge (50\Delta)!\ge 50!$, sufficiently large $n$ and $G,\bc$ as in Lemma~\ref{lem:linearly-balanced-vertices}. Then, property \emph{(c)} is satisfied with probability at least $0.8$.
\end{lemma}

\begin{proof}
Let the number of type-1 copies of $K_{2,2\Delta}$ in $G$ be $\beta n^{2\Delta + 2}$, where $\beta\in [1/(40\Delta)!, 1)$ by \Cref{lem:lots-of-gadget-sites}.
	Moreover, since there are at least $m/\Delta$ non-isolated vertices in $H$ and each of them has at most $\Delta^2+\Delta+1$ vertices at distance at most 2 in $H$, there must be a set $S$ of
	\[\frac{m/\Delta}{\Delta^2+\Delta+1}\cdot \frac{1}{\Delta}\ge 2\xi := 2\cdot \frac{m}{6\Delta^4}\]
	vertices with the same (positive) degree and pairwise distance at least 3 in $H$.
	We group these vertices into a set $M$ of $\xi$ pairs.
    Note that, for any copy of $H$ in $K_n$, the gadgets formed by the neighbourhoods of different pairs of vertices in $M$ are disjoint.
	
	Given a uniformly random copy $H'$ of $H$ in $K_n$ and a pair $(x,y)\in M$, denote by $I_{(x,y)}$ the indicator function of the event that the neighbourhoods of $x$ and $y$ form a positive gadget. Moreover, denote by $X$ the sum of the latter indicators.
	By a second-moment argument similar to~\eqref{eq:var}, we will show that, with probability at least $0.9$, we have $X\ge m/(45\Delta)!$ (the statement for the negative gadgets follows similarly).
	Indeed, let $p_1$ be the probability that two vertices $x,y$ with $(x,y)\in M$ have their neighbourhoods embedded into a type-1 copy of $K_{2,2\Delta}$ in a way that a positive gadget is formed, and let $p_2$ be the probability that two pairs in $M$ simultaneously satisfy this property. 
	Setting $\Pi := \prod_{i=0}^{2\Delta+1}(n-i)$, we obtain that 
	\[p_1 = \frac{\beta n^{2\Delta+2}}{\Pi} = \beta + O\biggl(\frac{\beta \Delta^2}{n}\biggr) \quad\text{and}\quad p_2 = p_1\cdot \frac{p_1\cdot \Pi - O(\Delta^2 n^{2\Delta+1})}{\Pi- O(\Delta^2 n^{2\Delta+1})} = p_1^2+O\bigg(\frac{\beta \Delta^2}{n}\bigg).\]
	Hence,
    \begin{align*}
        \Var(X) = \EE \bigg[ \bigg(\sum_{(x,y)\in M} I_{(x,y)}\bigg)^2 \bigg] - \xi^2 p_1^2 = \xi p_1 + \xi(\xi-1) p_2  -\xi^2 p_1^2 \leq \xi\beta + O\bigg(\frac{\xi^2 \beta \Delta^2}{n}\bigg),
    \end{align*}
	and Chebyshev's inequality shows that
	\[\mathbb P(X\le m/(45\Delta)!)\le \mathbb P(|X-\xi\beta|\ge \xi\beta/2)\le \frac{4\mathrm{Var}(X)}{\xi^2\beta^2}\le \frac{4}{\xi\beta} + O\bigg(\frac{\beta \Delta^2}{\beta^2 n}\bigg)\le 0.1,\]
	which shows the part of (c) concerning the positive gadgets.
	A similar argument for the negative gadgets and a union bound finishes the proof.
\end{proof}

We are ready to prove Theorem~\ref{thm:AlmostZeroSumDenseGraph}.

\begin{proof}[Proof of Theorem~\ref{thm:AlmostZeroSumDenseGraph}]
	Fix $m\ge (50\Delta)!$, $\nu\in (0, 1/(50\Delta)!]$, $\eps\le \nu^2/10^5$ and $\alpha = 1/5(45\Delta)!$.
    We begin by considering the case where the triplet $(H,G,\bc)$ is not inert, with the other case handled similarly.
    Let $b$ be the residue class of $m$ modulo $2$, and let $k \in [-\alpha m, \alpha m]$ such that $k\equiv b\,(\hspace{-0.62em}\mod 2)$.
    We then consider an embedding $H'$ of $H$ taken uniformly at random in $K_n$.
    As previously noted, property (a) holds with probability at least $0.8$ by Lemma~\ref{lem:balanced_colours} applied with $\delta = \nu/5\in (0,2^{-20}]$, and property (b) holds with probability at least $2/3$ by Markov's inequality.
    Additionally, property (c) is satisfied with probability at least $0.8$ by \Cref{lem:find-common-gadgets}.
    Therefore, properties (a), (b) and (c) hold simultaneously with positive probability, and we fix a copy $H_1$ of $H$ satisfying all of them.
    We proceed with four consecutive rounds of vertex exchanges to complete the~proof:
    \begin{enumerate}
        \item First, we select a type-2 copy $K$ of $K_{2,2\Delta}$ (as provided by Lemma~\ref{lem:2-gadget}) and perform up to $2\Delta+2$ exchanges of vertex pairs in $H_1$ to obtain a copy $H_2$ of $H$ in $K_n$ containing a 2-gadget embedded~into~$K$.
        We assume that the vertex exchange described in \Cref{rem:PositiveNegativeGadget} increases the value of $\bc(H_2)$ by $2$, the other case being analogous. 
        Moreover, combining the fact that $H_1$ and $H_2$ differ in at most $(2\Delta+2)2\Delta\le 8\Delta^2$ edges with property (b), it follows that $H_2$ contains at most $3\eps m+8\Delta^2$ non-edges of $G$.
        \item Next, by applying Lemma~\ref{lem:remove-ghost-edges} and Remark~\ref{rem:non-tree} (with $n\ge 5(\Delta+1)^3$, ensuring $n/2(\Delta+1)^2\ge 2\Delta +3$), we can find a copy $H_3$ of $H$ in $G$ that differs from $H_1$ in at most 
	    \[8\Delta^2 + 2\Delta(3\eps m+8\Delta^2)\le 24\Delta^3 + 6\Delta \eps m\]
	    edges. Furthermore, since the 2-gadget contains $2\Delta+2$ vertices, Remark~\ref{rem:non-tree} guarantees that $H_3$ can be obtained from $H_2$ while leaving this 2-gadget untouched. 
        \item Using property (c) and the inequality 
        \[\frac{m}{(45\Delta)!}\ge 5\bigg( (2\Delta+2)\Delta + 24\Delta^3 + 6 \Delta \eps m\bigg),\]
	it follows that there are at least $4m/5(45\Delta)!$ positive and at least $4m/5(45\Delta)!$ negative vertex-disjoint gadgets that remained untouched in the process of obtaining $H_3$ from $H_1$.
    Using \Cref{rem:PositiveNegativeGadget}, recall that the described vertex exchange in a positive (resp.\ negative) gadget changes the value of $\bc(H_3)$ by 2 or 4 (resp.\ $-2$ or~$-4$).
    Furthermore, we have that 
    \[\abs{\bc(H_3)} \leq 2\bigg(\frac{\nu m}{10}+(2\Delta+2)\Delta+ 24\Delta^3 + 6\Delta \eps m\bigg) \leq \frac{1}{2}\cdot \frac{4m}{5(45\Delta)!}\]
    Thus, $|k-\bc(H_3)| \leq 2m/5(45\Delta)!$. As a result, by activating a suitable number of positive or negative gadgets (depending on of the sign of $k-\bc(H_3)$), we obtain a copy $H_4$ of $H$ in $G$ such that $k-\bc(H_4)=0$ or $k-c(H_4)=2$.
        \item If $k-\bc(H_4)=0$, we are done. If $k-\bc(H_4)=2$, then we perform a final vertex exchange as described in \Cref{rem:PositiveNegativeGadget} within the 2-gadget embedded into $K$ during the first round of vertex exchanges. 
        This way, we obtain a copy $H_5$ of $H$ in $G$ satisfying $\bc(H_5)=k$, as desired.\qedhere
    \end{enumerate}
\end{proof}

\section{Concluding remarks and open problems}\label{sec:conc}

We begin with a few comments concerning analogues of our results for digraphs.
The topic was also raised by Freschi, Hyde, Lada, and Treglown~\cite{FHLT21} who asked whether any $n$-vertex $r$-edge-coloured digraph with minimum in- and out-degree at least $(1/2+1/2r+o(1))n$ contains a Hamilton cycle with significant colour-bias.\footnote{Recall that, in a digraph containing vertices $u$ and $v$, the directed edges $uv$ and $vu$ can appear simultaneously.} 
The next simple example shows that, in fact, no in- and out-degree conditions are sufficient to guarantee a significantly colour-biased copy of any rooted spanning tree oriented away from its root. 
Since Hamilton cycles are one edge away from a consistently oriented spanning path, the example also gives a negative answer to the question in~\cite{FHLT21}.

\begin{example}
Partition a set $V$ of $n$ vertices into $r$ parts $V_1, \ldots, V_r$ of sizes as equal as possible. In the complete digraph on $V$, 
colour every edge towards $V_i$ in colour $i$.
For every rooted spanning tree $T$ with edges oriented away from its root, every vertex except the root has a single incoming edge.
As a result, every copy of $T$ has between $\lceil n/r\rceil$ and $\lfloor n/r\rfloor - 1$ edges in each colour.
\end{example}

We note that the above example relies on the fact that almost all vertices in the tree $T$ have the same in-degree.
Going further by analysing spanning trees and, in general, bounded degree spanning subgraphs with in- and out-degree sequences which are far from regular remains an interesting open question.

We believe that some aspects of our work can be pushed further. 
First or all, we observe that the polynomial dependency on $r$ in \Cref{col:discErdos} is necessary.
To see this, fix $k$ satisfying $k^2+k\le 2r$ and consider a copy of $K_n$ with vertex partition $V_1, \dots, V_k$ of sizes as equal as possible.
Colouring the edges between $V_i$ and $V_j$ (with $i,j$ possibly coinciding) in a distinct colour for all pairs $i,j$ provides a construction where no spanning tree has discrepancy more than $2\lceil n/k\rceil$.
The exact dependency on $r$ remains a curious open problem.
Regarding \Cref{thm:col_disc_3}, we believe that ideas of Komlos, Sarkozy and Szemerédi~\cite{komlos1995proof} can be used to remove the dependency of $\eps$ on $\Delta$ but no such attempt has been made in our work.

Another intriguing direction of research is to understand whether \Cref{thm:AlmostZeroSumDenseGraph} could be generalised for $r\ge 2$ colours.
While we believe that extending our method to the 3-coloured setting is possible, the main approach breaks down for $r\ge 4$ due to pathological constructions similar to \Cref{ex:1}.
In this direction, Pardey and Rautenbach~\cite{pardey2022almost} conjectured that in a balanced colouring of the complete graph with many colours, one can always find a perfect matching which is almost colour balanced. 
This conjecture is the focus of a paper by the first author~\cite{Hol24}.
However, understanding the general many-colour setting when $H$ is a bounded-degree graph, or even the simpler case of an $R$-factor for some fixed graph $R$, remains widely open.

\paragraph{Acknowledgements.} The first and third authors were funded by Trinity College, Cambridge.
The second author was supported by the Austrian Science Fund (FWF) via ESPRIT grant No. 10.55776/ESP624.
The fourth author was funded by EPSRC (Engineering and Physical Sciences Research Council) and by the Cambridge Commonwealth, European and International Trust.
The second and fourth authors acknowledge the support of the Bulgarian Ministry of Education and Science, Scientific Programme ``Enhancing the Research Capacity in Mathematical Sciences (PIKOM)'', No. DO1-241/15.08.2023, and are grateful to the Institute of Mathematics and Informatics at the Bulgarian Academy of Sciences for their hospitality.

\bibliographystyle{abbrvnat}  
\renewcommand{\bibname}{Bibliography}
\bibliography{bib}
	
\end{document}